\documentstyle[12pt]{article}
\newcommand{\real}{\mbox{I$\!$R}}

\newcommand{\al}{\alpha}

\newcommand{\de}{\delta}
\newcommand{\th}{\theta}

\newcommand{\ep}{\epsilon}
\newcommand{\qed}{\mbox{\bf \ \ \ \ \ Q.E.D.}\\[0.2in]\hspace*{0.5in}}
\newcommand{\bk}{\\[0.01in] \hspace*{0.5in} }
\begin{document}

\title{Rotationally Symmetric $p$-harmonic maps}
\author{Man Chun Leung\\Department of Mathematics, National University of
Singapore,\\ Singapore 119260  \ \ {\it matlmc@nus.sg}}
\date{April 1996}
\maketitle
\begin{abstract}This is a report on our study of a second order ordinary
differential equation corresponding to rotationally symmetric 
$p$-harmonic maps.
We show unique continuation and Liouville's type theorems for positive 
solutions.
Asymptotic properties of the positive solutions are investigated. 
\end{abstract}
\vspace{0.2in} KEY WORDS: $p$-harmonic map, positive solution,  rotationally
symmetric\\[0.1in]   1991 AMS MS Classifications: 34C11, 34D05, 58E20

\vspace{0.5in}

{\bf \Large 1. \ \ Introduction}

\vspace{0.3in}

In this report we study $p$-harmonic maps between {\it model} Riemannian
manifolds. For $n \ge 2$, let $(S^{n - 1}, d\vartheta^2)$ be the unit 
sphere in
${\real}^n$ with the induced Riemannian metric $d\vartheta^2$. Let $f$ 
and $g$
be functions in $C^2 ([0,
\infty))$ which satisfy the following
$$
f (0) = g (0) = 0\,, \ \ \ \ f' (0) = g' (0) = 1 \ \ \
{\mbox{and}}
\ \ \  f (r), g(r) > 0\,. \leqno (1.1)
$$ 
Consider the
following  model Riemannian manifolds [5] 
$$ M (f) = ( [0, \infty ) \times S^{n-1}\,, \ dr^2 + f^2(r) d\vartheta^2 
\ )\,,
$$      
$$ 
N (g)=( [0, \infty ) \times S^{n-1}\,, \ dr^2 + g^2(r)d\vartheta^2\ )\,,
$$ 
where $f, g \in C^2 ([0, \infty))$ satisfy the conditions in (1.3). The
Euclidean space and the hyperbolic space are corresponding to $f (r) = r$ 
and $g
(r) = \sinh r$, respectively. A map $F: M^n (f) \to M^n (g)$ is called a
rotationally symmetric map if 
$$F (r, \vartheta) = (\alpha (r)\,, \vartheta) \ \ \ \ {\mbox{for \ \ 
all}} \ \ \
r > 0 \ \ \ {\mbox{and}} \ \ \vartheta \in S^{n - 1}\,,$$ 
where $\alpha : [0, \infty) \to [0, \infty)$ with $\alpha (0) = 0$.  The 
energy
density [4] of a rotationally symmmetric map $F: M^n (f) \to M^n (g)$ is 
given by 
$$
\theta (r) = (\alpha' (r))^2 + (n - 1) {{g^2 (\alpha (r))}\over {f^2 (r)}}\,.
\leqno (1.2)
$$   
Given $p > 0$, the $p$-energy functional [13,\,14] for
$C^2$ rotationally symmetric maps is given by 
$$
\int_0^\infty \left\{ \alpha'^2 (r) + (n - 1) {{g^2 (\alpha (r))}\over {f^2
(r)}} \right\}^{p\over 2} f^{n -1} (r) dr\,. \leqno (1.3)
$$ 
It can be shown that equation the Euler-Lagrange equation of the
$p$-energy functional is 
\begin{eqnarray*} (1.4) \ \ \ \ \  \ \ \theta^{{p\over 2} - 1} (r) 
\alpha''(r) &
+ & [(n - 1)
\theta^{{p\over 2} - 1}(r) {{f'(r)}\over {f(r)}} +  (\theta^{{p\over 2} -
1})'(r) ] \alpha'(r)\\ & \ & \ \ \ \ \ \ \ \ \ \ \ \ \ \ -  (n - 1)
\theta^{{p\over 2} - 1} (r){{ g (\alpha (r)) g' (\alpha (r))}\over {f^2 
(r)}} =
0\,.\ \ \ \ \ \ \ \ \ \ \ \ \  \ \ 
\end{eqnarray*}
Critical points of the $p$-energy functional (1.4) are known as rotationally
symmetric $p$-harmonic map from $M (f)$ to $N (g)$. We remark that the 
notion 
of $p$-energy, which generalizes the usual energy as in harmonic maps 
[4], was
first  introduced by K. Uhlenbeck [13]. The existence question has been 
studied
later by various authors, using the  direct method of Calculus of 
Variations. In
[14], White studied maps  minimizing the $p$-energy in certain homotopy 
classes.
Xu and Yang  [15] considered the existence of $p$-harmonic maps between 
spheres, 
generalizing the works of previous authors on harmonic maps between  
spheres. It
is worth mentioning that the direct method yields normally  only finite 
energy
maps. For noncompact manifolds, Nakauchi and Takakuwa [10] study 
concentration
compactness and gap theorems for $p$-harmonic maps. \bk
For $p = 2$, equation (1.4) becomes
$$\alpha''(r) + (n - 1) {{f'(\alpha)}\over {f(\alpha)}} \alpha' (r)  -  
(n -
1){{ g (\alpha (r)) g' (\alpha (r))}\over {f^2 (\alpha)}} = 0\,. \leqno 
(1.5)$$
Equation (1.5) corresponds to rotationally symmetric harmonic maps, which 
have
been studied by Ratto and Rigoli [11] and Cheung and Law [1].
In this note we discuss the existence and properties of positive solutions
$\alpha \in C^2 ( (0, \infty)$ to equation (1.4) with
$$\lim_{r \to 0^+} \alpha (r) = 0\,.$$
The results in this note are contained in the works of Cheung-Law-Leung 
[2] and
also in [8,\,9].\bk
We first show that any positive solution must be increasing. Equation 
(1.4) is a
quasilinear second order differential equation, where the nonlinear term
$\theta^{p/2 - 1} (r)$ makes the equation complicated. We introduce a 
first order
equation on $\theta (r)$ and $\alpha (r)\,$ [equation (2.3)]. From this 
we study
unique continuation of nonnegative solutions to equation (1.4). For
$p > 2$, we show that if  $\alpha (r)
\in C^2 (0, \epsilon)$ is a nonnegative solution to (1.4) with $\lim_{r 
\to 0^+}
\alpha (r) = 0$, and $\alpha (r)$  goes to zero at arbitrary high order, then
$\alpha
\equiv 0$ on $(0, \epsilon)$ (lemma 2.6). This unique continuation 
property is
well-known for harmonic maps [7] but is unknown to us for $p$-harmonic 
maps in
general with $p \not=2$ (c.f. [6]). By refining the argument, we can estimate
the order that
$\alpha$ can approach zero. We obtain a Liouville's type theorem for certain
rotationally symmetric
$p$-harmonic maps into the Euclidean space or the hyperbolic space (theorem
2.17 \& 2.23).\bk
The differential equation (1.4) is solved locally near zero by solving   an
initial value problem for small $r$ (theorem 3.1). This is achieved by
rewriting the differential equation (1.4) as a system of two integral  
equations
and solving it by comparing it to the Euler equation. This local solution is
shown to be extendible to the entire positive real axis under certain
conditions on $f$ and $g$. We study asymptotic properties of rotationally
symmetric $p$-harmonic maps from the hyperbolic space to itself, that is, $f
(r) = g (r) = \sinh r$. The identity map Id $: M^n (f) \to M^n (f)$ 
corresponding
to
$\alpha (r) \equiv r$ is a rotationally symmetric $p$-harmonic map. We
investigate asymptotic properties of a rotationally symmetric 
$p$-harmonic map
by seeing whether there is a point "above" or "below" the identity map. More
precisely, for $p > 2$, we show that if $F: M^n (\sinh r) \to M^n (\sinh 
r)$ is a
rotationally symmetric $p$-harmonic map from the hyperbolic space to 
itself and
if there is a point
$r_o > (\ln 3)/2$ such that $\alpha (r_o) < r_o$ and $\alpha' (r_o) < 1$, 
then
$\alpha$ is a bounded function. On the other hand if $\alpha (r_o) > r_o$ and
$\alpha' (r_o) > 1$, then for all positive number C, $\alpha (r) > r + C$ for
all $r$ large. If there is a positive constant $C$ such that $r - C \le 
\alpha
(r) \le r + C$ for all $r$ large, then $\alpha$ is asymptotic to the line
$y = x$ (theorem 4.18). \bk
We consider the asymptotic properties of positive solutions to (1.4) when 
$M (f)$
is similar to the hyperbolic space and
$N (g)$ is similar to the Euclidean space. we show that $\alpha (r)$ 
either goes
to infinity exponentially or it is bounded (theorem 5.12). If $g$ grows 
at most
exponentially, then we show that either $\alpha$ grows at least linearly or
$\alpha$ is bounded (theorem 5.19). We show that, under certain 
conditions on $f$
and
$g$, all rotationally symmetric
$p$-harmonic maps with $p > 2$ are bounded (corollary 3.24).\bk

\vspace{0.5in}

{\bf \Large 2. \ \ Unique continuation and Liouville's type theorem}

\vspace{0.3in}

We begin this section with local considerations of equation 
(1.1).\\[0.1in] {\bf
Lemma 2.1. [2]} \ \ {\it Let $y (r) \in C^2 (0, R)$ be a positive 
solution to
(1.1) with $\lim_{r \to 0^+} y (r) = 0$, where
$R$ is a positive number. Then
$y' (r) > 0$ for all $r \in (0, R)$.}\\[0.2in]   {\bf Proof.} \ \ Since 
$y (r) >
0$ for all $r > 0$ and $\lim_{r \to 0^+} y (r) = 0$, we can find a point
$r_o$ such that $y' (r_o) > 0$. Suppose that there exists a point $r_1 > r_o$
such that $y' (r_1) \le 0$, then there exists a point $r' \in (r_o, r_1]$ 
such
that $y (r') > 0$, $y' (r') = 0$ and $y'' (r') \le 0$. At $r'$, equation 
(1.1)
gives
$$\theta^{{p\over 2} - 1} (r') y''(r') = (n - 1)
\theta^{{p\over 2} - 1} (r'){{ g (y (r')) g' (y (r'))}\over {f^2 
(r')}}\,,$$ and 
$$\theta (r') = (n - 1) {{g^2 (y (r'))}\over {f^2 (r')}} + (y' (r'))^2 > 
0\,.$$
Thus $y'' (r') > 0$, contradiction. Therefore $y' (r) > 0$ for all $r \in
(r_o\,, R)$. The proof is completed by noticing that we can let $r_o
\to 0$.\qed
For $p > 0$, let 
$$\Theta (r) = \theta^{{p\over 2} - 1} (r)\,, \ \ \ \ r > 0\,.\leqno 
(2.2)$$ Then
$$\Theta' (r) \theta = ({p\over 2} - 1) \Theta (r) [(n - 1) ({{g^2 (\alpha
(r))}\over {f^2 (r)}})' + 2 \alpha' (r) \alpha'' (r)]\,.$$ Using equation 
(1.1)
to replace the term $\Theta (r) \alpha '' (r)$ we have 
\begin{eqnarray*}(2.3)& \ &\Theta' (r) [ ( p - 1 ) (\alpha' (r))^2 + (n - 1)
{{g^2 (\alpha (r) )}\over {f^2 (r)}} ]\\ & = & (p - 2) (n - 1)
\Theta (r) \{ {{2 g (\alpha (r) ) g' (\alpha (r))
\alpha' (r)}\over {f^2 (r)}} - {{f' (r)}\over {f (r)}} [ {{ g^2 (\alpha
(r))}\over {f^2 (r)}} + (\alpha' (r))^2 ] \}\,.
\end{eqnarray*} Equation (2.3) involves only first derivatives of 
$\Theta$ and
$\alpha$. For some $\epsilon > 0$, assume that a nonnegative $C^2$-solution
$\alpha$ of (1.4) exists on
$(0, \epsilon)$ with $\lim_{r \to 0^+} \alpha (r) = 0$. Given $\kappa > 
0$, we
can find positive constants
$a_o, b_o$ and
$c_o$ depending on $\kappa$ such that
$$
|f' (r) | \le a_o\,, \ \ \ f (r) \ge c_or \ \ \ {\mbox{and}} \ \ \ |g' (y 
(r))|
< b_o \ \ \ \ {\mbox{for \ \ all}} \ \ \ r \in (0, {{\kappa\over {\kappa 
+ 1}}
\epsilon})\,.
$$ 
Then (2.3) gives
\begin{eqnarray*} (2.4) \ \ \ \ \ \ \ \ \ \ \ \ \ \ \ & \ & \ \Theta' (r) 
[ ( p -
1 ) (\alpha' (r))^2 + (n - 1) {{g^2 (\alpha (r) )}\over {f^2 (r)}} ]\\  & 
\le &
(p - 2) (n - 1) \Theta (r) {{(a_o + b_o)}\over {c_or}}[ {{ g^2 (\alpha
(r))}\over {f^2 (r)}} + (\alpha' (r))^2 ]\,.\ \ \ \ \ \ \ \ \ \ \ \ \ 
\end{eqnarray*} Therefore there is a positive constant $C = C
(a_o,b_o,c_o,p,n,\kappa)$ such that for
$p > 2$ we have
$$
\Theta' (r) \le {C \over r} \Theta (r) \ \ \ \ \ \ {\mbox{on}} \ \ \ (0,
{{\kappa\over {\kappa + 1}} \epsilon})\,.\leqno (2.5)
$$ We have the following
unique continuation property for solutions of (1.1).\\[0.2in]  
{\bf Lemma 2.6. [9]}
\ \ {\it For $p > 2$, let $\alpha (r) \in C^2 (0, \epsilon)$ be a nonnegative
solution to (1.1) with $\lim_{r \to 0^+} \alpha (r) = 0$. If $\alpha (r) 
= O
(r^k)$ near $0$ for all $k > 0$, then $\alpha
\equiv 0$ on
$(0, \epsilon)$}\\[0.1in] {\bf Proof.} \ \ Assume that $\alpha \not\equiv 
0$ on
$(0, \epsilon /2)$. We first show that $\alpha$ cannot be zero on $(0, 
\delta)$
for any $\delta \in (0,
\epsilon /2)$. Suppose that $\alpha \equiv 0$ on
$(0, \delta)$ for some $\delta \in (0, \epsilon /2)$. Since $\alpha 
\not\equiv 0$
on $(0, \epsilon /2 )$, we may assume that $\alpha (r) > 0$ on $(\delta,
\epsilon /2)$. Integrating (2.5) we have
$$\ln \Theta (r) |^a_{b} \le C \ln r |^a_b\,,$$ where $\delta < b < a < 
\epsilon
/2$. That is 
$$ 
\Theta (a) \le \Theta (b) ({a\over b})^C\,.\leqno (2.7)
$$ Let $b \to \delta >
0$, we have $\Theta (b) \to 0$, but $\Theta (a) > 0$, contradicting 
(2.7). Thus
$\alpha$ cannot be zero on $(0, \delta)$ for any
$\delta \in (0, \epsilon /2)$. Thus we can find a point
$r_o \in (0, \epsilon /2)$ such that $\alpha (r_o) > 0$ and $\alpha' 
(r_o) > 0$.
Suppose that there exists a point $r_1 \in (0, \epsilon /2)$
$r_1 > r_o$ such that $\alpha (r_1) = 0$, then there exists a point $r' \in
(r_o, r_1)$ such that $\alpha (r') > 0$, $\alpha' (r') = 0$ and $\alpha'' (r)
\le 0$. At $r'$, equation (1.4) gives
$$\theta^{{p\over 2} - 1} (r') \alpha''(r') = (n - 1)
\theta^{{p\over 2} - 1} (r'){{ g (\alpha (r')) g' (\alpha (r'))}\over {f^2
(r)}}\,,$$ and 
$$\theta (r') = (n - 1) {{g^2 (\alpha (r'))}\over {f^2 (r')}} + (\alpha' 
(r'))^2
> 0\,.$$ Therefore we have $\alpha'' (r') > 0$, contradiction.  Therefore
$\alpha (r) > 0$ for all $r \in (r_o, \epsilon /2)$. As $\lim_{r \to 0^+} 
\alpha
(r) = 0$ and
$\alpha
 \not\equiv 0$ on $(0, \delta)$ for all $\delta > 0$, we can let
$r_o \to 0$. Thus $\alpha (r) > 0$ on $(0,
\epsilon /2)$. Similarly we can show that $\alpha' (r) > 0$ on $(0, \epsilon
/2)$. Integrating (2.5) we obtain 
$$ \Theta (a) \le \Theta (r) ({a\over r})^C\,,$$ where in this case 
$\epsilon /2
> a > r > 0$ and $a$ is a constant. We have
$$
\Theta (r) \ge C' r^C \leqno (2.8)
$$ 
for some positive constant $C'$ That is,
$$
\theta (r) = (n - 1) {{g^2 (\alpha (r))}\over {f^2 (r)}} + (\alpha' 
(r))^2 \ge
C'' r^{{2C}\over {p - 2}}\,,\leqno (2.9)
$$ 
where $C''$ is a positive constant. Take a number $k$ such that 
$$
k > {{C}\over {p - 2}} + 1\,.
$$  
Since $\alpha (r) \le C_k r^k$ for $r$ small,
where
$C_k$ is a positive constant, and $f (r), g(r) \sim r$ when $r$ is small, 
we have
$$ (n - 1) {{g^2 (\alpha (r))}\over {f^2 (r)}} \le C''' r^{2k - 2}\,,
$$ where $C'''$ is a positive constant.  Hence by (2.8) we have
$$ (\alpha' (r))^2 \ge C_o r^{{2C}\over {p - 2}} \leqno (2.10)
$$ 
for some positive constant $C_o$. Therefore we have
$$
\alpha (r) - \alpha (b) = \int^r_b \alpha' (s) ds \ge C'' \int_b^r
s^{{2C}\over {p - 2}} ds = C_o (r^{ {{2C}\over {p - 2}} +1 } - b^{ 
{{2C}\over {p
- 2}} + 1})\,, \ \ \ {\epsilon \over 2} > r > b > 0\,.
$$ 
As $\lim_{b \to 0^+} \alpha (b) = 0$, we have
$$
\alpha (r) \ge C_1 r^{ {{2C}\over {p - 2}} +1 }\,, \leqno (2.11)
$$
contradicting that $\alpha (r) \sim O (r^k)$ for all $k > 0$. Therefore 
we have
$\alpha (r) \equiv 0$ on $(0, \epsilon /2)$. Similar argument shows that 
$\alpha
(r)
\equiv 0$ on $0, ({\kappa\over {\kappa + 1}} \epsilon )$ for all $\kappa 
> 0$.
Let $\kappa \to \infty$, we have $\alpha (r) \equiv 0$ on $(0,
\epsilon)$.\qed   Let $\alpha \in C^2 ((0, \epsilon))$ be a nonnegative 
solution
to $(1.1)$ with
$\lim_{r \to 0^+} \alpha (r) = 0$. If we assume that $\alpha \not\equiv 
0$, then
the proof of lemma 2.6 shows that $\alpha (r) > 0$, hence $\alpha' (r) > 
0$ on
$(0,
\epsilon)$ [2]. Given any $\delta > 0$, as $\lim_{r \to 0^+} \alpha (r) = 
0$, we
can find $\epsilon_o < \epsilon$ such that 
$$0 < \alpha (r) < \delta \ \ \ \ {\mbox {on}} \ \ \ \ (0, 
\epsilon_o)\,.$$ As
$f, g \in C^2 (0, \infty)$ with $f' (0) = g' (0) = 1$, given
$\varepsilon_1 \in (0, 1)$, we can assume that
$\delta$ and $\epsilon_o$ are small enough such that 
$$0 < f' (r)  \le 1 + \varepsilon_1\,, \ \ \ \  f (r) \ge (1 - 
\varepsilon_1) r
\ \ \ {\mbox{and}} \ \ \ |g' (y (r))| < 1 + \varepsilon_1 \ \ \ \ 
{\mbox{for \ \
all}} \ \ \ r \in (0, \epsilon_o)\,.$$ If we assume that $p > 2$ and $n 
\ge 2$,
then (2.2) gives
$$
\Theta' (r) \le (p - 2) (n - 1) ({{1 + \varepsilon_1}\over {1 -
\varepsilon_1}} ) \Theta (r) \ \ \ \ {\mbox {on}} \ \ \ \ (0, \epsilon_o)\,.
\leqno (2.12)
$$  That is, $C = (p - 2) (n - 1) (1 + \varepsilon_1)/(1 -
\varepsilon_1)$ in (2.4). It follows as in the proof of lemma 2.5 and 
inequality
(2.9) that 
$$
\alpha (r) \ge C_1 r^{ {{2C}\over {p - 2}} +1 } = C_1 r^{ 2 (n - 1) ({{1 +
\varepsilon_1}\over {1 - \varepsilon_1}}) + 1 }\,. \leqno (2.13)
$$ As $\varepsilon_1$ can be made arbitrarily small, therefore we have the
following refinement of lemma 2.6.\\[0.2in]  
{\bf Theorem 2.14 [9]} \ \ {\it For $p >
2$ and $n
\ge 2$, let $\alpha (r) \in C^2 (0,
\epsilon)$ be a nonnegative solution to (1.1) with $\lim_{r \to 0^+} 
\alpha (r)
= 0$. If there exist constants $k > 2n - 1, C' > 0$ and $\epsilon_o \in (0,
\epsilon)$ such that , 
$\alpha (r) \le C' r^k$ for all $r \in (0, \epsilon_o)$, then
$\alpha \equiv 0$ on $(0, \epsilon)$}\\[0.1in]
{\bf Lemma 2.15. [2]} \ \ {\it For $p > 2$ and $n \ge 2$, let $\alpha (r) 
\in C^2
(0, R)$ be a positive solution to (1.1) with $\lim_{r \to 0^+} \alpha (r) 
= 0$,
where
$R$ is a positive number. Suppose that $f$ and $g$ satisfy the conditions in
(1.3) and $|f' (r)| \le a$ and
$|g' (\alpha (r))| \leq b$ for all $r \in (0, R)$, where $a, b > 0$ are 
constants,
then the relative energy density $\theta (r)$  satisfies the following
differential inequality on $(0, R)$:}
$$-\frac{(n-1)(p - 2)(a + b)}{(p - 1)}\, {1\over {f (r)}} \leq
\frac{{(\th^{q-1})}' (r)}{\th^{q-1} (r)}
\leq  {{(n-1)(p - 2)(a + b)}\over {\min \,\{ p - 1\,, \ n - 1 \} }} \, 
{1\over
{f (r)}}\,. \leqno (2.3)$$  {\bf Proof.} \ \ By (2.3)  we have 
\begin{eqnarray*}  & \ & \ (\theta^{q - 1})' (r) [ ( p - 1 ) (\alpha'
(r))^2 + (n - 1) {{g^2 (\alpha (r) )}\over {f^2 (r)}} ]\\ & = & (p - 2) 
(n - 1)
\theta^{q - 1} (r) \{ {{2 g (\alpha (r) ) g' (\alpha (r)) \alpha' 
(r)}\over {f^2 (r)}} - {{f'
(r)}\over {f (r)}} [ {{ g^2 (\alpha (r))}\over {f^2 (r)}} + (\alpha' 
(r))^2 ] \}\,.
\end{eqnarray*}  
As $\alpha'> 0$, we have
\begin{eqnarray*} (2.16) \ \ \ \ \ \ \ \ \ \ & \ & {{2 g (\alpha (r) ) g' 
(\alpha
(r))
\alpha' (r)}\over {f^2 (r)}} - {{f' (r)}\over {f (r)}} [ {{ g^2 (\alpha
(r))}\over {f^2 (r)}} + (\alpha' (r))^2 ] \\ & \le & {1\over {f (r)}} \{ b
[{{g^2 (\alpha (r) )}\over {f^2 (r)}} + (\alpha' (r))^2] + a [{{g^2 
(\alpha (r)
)}\over {f^2 (r)}} + (\alpha' (r))^2] \} \\  & \le & {1\over {f (r)}} (a 
+ b) 
[{{g^2 (\alpha (r) )}\over {f^2 (r)}} + (\alpha' (r))^2]\,.  \ \ \ \ \ \ 
\ \ \ \ 
\ \ \ \ \ \ \ \ \ \ \ \ \ \ \ \ \ \ \ \ \ \ \ \ \ \ \ \ \ \ \ 
\end{eqnarray*} Therefore we have
$$(\theta^{q - 1})' (r) \le C' {{\theta^{q - 1} (r)}\over {f (r)}}$$ for 
$r \in 
(0, R)$, where $C = (p - 2) (n - 1) (a + b)/ \min \{ p - 1\,, \ n - 1 \}$ 
is a
positive constant. On the other hand, by (2.4) we have
\begin{eqnarray*} (p - 1) \th {(\th^{q-1})}'&\geq&-\frac{(n-1)(p-2) (a +
b)}{f(r)}
\, \th^{q-1} [{(\alpha')}^2 + \frac{g(\alpha)^2}{f(r)^2}]\, \ ,\\
 &\geq&-\frac{(n-1)(p - 2)(a + b)}{f(r)}\, \th^q\ .
\end{eqnarray*} Thus
$$ {(\th^{q-1})}'\geq -\frac{2(n-1)(q-1)(a + b)}{(p - 1)f(r)}\, 
\th^{q-1}\ .
$$
{\bf Q.E.D.}\bk
For $p > 2$ and $n \ge 2$, let $F=(\phi , \alpha(r))$ be a rotationally 
symmetric
$p$-harmonic map from a model space $M^n(f)$ to
$N^n (g)$, with $\alpha (r) \ge 0$ defined on $(0, \infty)$ and $\lim_{r 
\to 0^+} \alpha
(r) = 0$, where $f$ and $g$ satisfy the condition in (1.3). We have the 
following
Liouville's type theorem for $p$-harmonic maps.\\[0.1in]  
{\bf Theorem 2.17. [2]} \ \
{\it Assume that $|f' (r)| \le c_o^2,$ $g' (r) > 0$ and $g'^2 (r) - g (r) 
g'' (r)
\ge c^2$ for all $r > 0$, where
$c_o$ and $c$ are positive numbers. If $\alpha$ is bounded, then $\alpha 
\equiv
0$.}\\[0.1in]
{\bf Remark.} \ \ In
the Euclidean and hyperbolic case, we have $g (r) = r$ and
$g (r) = \sinh r$, respectively. In both cases, we have $g'^2 (r) - g (r) g''
(r) = 1$ and hence the condition in the above theorem is satisfied. The
condition on $f$ is satisfied for the type of functions $f (r) = r^\alpha +
\sin r$ for $r \gg 1$ and $\alpha \in (0, 1]$ (c.f. [10] for Liouville's type
theorems for harmonic maps).\\[0.1in]  {\bf Proof.} \ \ Assume that 
$\alpha$ is
bounded and
$\alpha
\not\equiv 0$, then lemma the proof of lemma 2.6 shows that $\alpha (r) > 0$,
$\alpha' (r) > 0$ and hence
$\theta (r) > 0$ for all
$r > 0$. The
$p$-harmonic map equation is given by
$$
(\theta^{q-1} )'\alpha' +\theta^{q-1}  \alpha'' 
+(n-1)\frac{f'}{f}\theta^{q-1}
\alpha'-(n-1)\frac{\theta^{q-1}}{f^2} g(\alpha) g' (\alpha) =0\,.
$$ 
Since $y$ is bounded, we
may assume that 
$|g' (\alpha (r))| < a$ for some constant $a > 0$. As in (2.16) we have
$$(\theta^{q - 1})' (r) < C {{\theta^{q - 1} (r)}\over {f (r)}}$$ for $r 
> 0$,
where $C = (p - 2) (n - 1) (a + c_o^2)/ \max \,\{ p - 1\,, \ n - 1 \} + 
1$ is a
positive constant. We have
$$ C \frac{\theta^{q-1}}{f}\alpha'+\theta^{q-1} \alpha'' +
(n-1)\frac{f'}{f}\theta^{q-1}\alpha'-(n-1)\frac{\theta^{q-1}}{f^2} 
g(\alpha) g' (\alpha)> 0\ .
$$ Since $\theta \ne 0$ except at the origin, this simplifies to
$$ {{\alpha''}\over {g(\alpha) g' (\alpha)}}  + \frac{n-1}{f}(C +f')\, 
{{\alpha'}\over { g(\alpha) g'
(\alpha)}} -(n-1)\frac{1}{f^2} > 0\,, \leqno (2.18)
$$ 
for all $r \ge 1$. For $r \ge 1$, let $H = \alpha'/\,[g (\alpha) g' 
(\alpha)]$.
We have 
$$ 
{{\alpha''}\over {g(\alpha) g' (\alpha)}} = H' + H^2 [g'^2 (\alpha) + g 
(\alpha)
g'' (\alpha)]\,.
$$ 
Since $\alpha$ is a bounded positive solution, we can find a positive 
constant
$C_o > 1$ such that $C_o^2 \ge g'^2 (\alpha) + g (\alpha) g'' (\alpha)$ 
for all $r \ge 1$. Thus
(2.18) gives
$$ C_o^2 H^2 + H' +\frac{n-1}{f}(C +f')\, H -\frac{(n-1)}{f^2} > 0.
\leqno (2.19)
$$ Consider the quadratic form
$$ C_o^2 H^2 + \frac{n-1}{f}(C +f') H -\frac{(n-1)}{f^2}\ .
$$ It is nonpositive for $H(r)\in [0,b_C(r)]$ where $b_C(r)$ is given by
$$  b_C (r)= {{n-1}\over {2fC_o^2}} \{ [ (C +f')^2 + {{4C_o^2}\over 
{n-1}}]^{1/2}
-(C+f') \} > 0\,. \leqno (2.20)
$$ Hence
$$ H'(r) > 0\quad\mbox{ whenever}\ \quad H(r)\leq b_C(r)\ .
\leqno (2.21)
$$ Consider the function $q (x) = \sqrt{x^2 + b^2} - x > 0$, where $b$ is a
positive number. If $x \le 0$, then $q (x) \ge b$. For $x > 0$, we have 
$q' (x)
< 0$. Thus $q (x) \ge c_1^2$ if $|x|$ is bounded, where $c_1$ is a positive
number. As $|f'| \le c_o^2$, we have $|C + f'| \le C + c_o^2$ and $f (r) 
\le C_1
r$ for some positive constant $C_1 > 0$, we can find a positive number
$C_2$ such that   
$$b_C (r) \ge {{C_2}\over r}\,.$$ Thus 
$$\int_1^\infty b_C (r) dr = \infty\,.$$ On the other hand we have 
$$ - ( \ln {{g' (\alpha)}\over {g (\alpha)}} )' = H [g'^2 (\alpha) - g 
(\alpha) g'' (y)] \ge  c^2 H
$$ for $r \ge 1$. Thus 
$$ - \ln {{g' (\alpha (r))}\over {g (\alpha (r))}} |^\infty_1 = - 
\int_1^\infty ( \ln {{g'
(\alpha)}\over {g (\alpha)}} )' dr \ge c^2 \int_1^\infty H (r) dr\,.
$$ As $g' (r) > 0$ for all $r > 0$ and $y$ is bounded and positive for $r \ge
1$, we have 
$$
\int_1^{\infty} H(r) dr < \infty\ .
$$ Since $H \geq 0$ and $\int_1^{\infty} b_C (r)dr =\infty$ , for some
sufficiently large $r$ we have $H(r)\leq b_C(r)$ and $H'(r)\leq 0$ (c.f. 
[6]),
contradicting (2.21).\qed
It follows from (2.3) that if $f' (r) \ge 0$, then we have
$$
\Theta' (r) [ ( p - 1 ) (\alpha' (r))^2 + (n - 1) {{g^2 (\alpha (r) 
)}\over {f^2
(r)}} ] \le (p - 2) (n - 1)
\Theta (r) {{2 g (\alpha (r) ) g' (\alpha (r))
\alpha' (r)}\over {f^2 (r)}}\,.
$$ In addition, if $p > 2$ and and $g' (r) \le a g (r)$ for all
$r > 0$, where $a$ is a positive constant, then there exists a constant $C_o$
such that  
$$
\Theta' (r) \le C_o \Theta (r) \alpha' (r)\,. \leqno (2.22)
$$ 
We use this to
proof the following Liouville's type theorem.\\[0.2in] 
{\bf Theorem 2.23 [9]} \ \
{\it Assume that there exists positive constants $c_o$ and $c_1$ such 
that $0
\le f' (r) \le c_o^2$ and
$f (r)
\le c_1$ for all $r \ge r_o$. For
$p > 2$, let
$F (r,
\vartheta) = (\alpha (r)\,,
\vartheta)$ is a rotationally symmetric $p$-harmonic maps from $M (f)$ to the
Euclidean space or the hyperboic space. If there exists a positive 
constant $C'$
such that $\alpha' (r)
\le C'$ for all
$r > 0$, then $\alpha \equiv 0$.}\\[0.1in] 
{\bf Remarks.} \ \ In fact we prove a
more general statement. The conditions below  on $g$ are satisfied for the
Euclidean and hyperbolic metrics, where $g (r) = r$ and $g (r) = \sinh r\,,$
respectively.  Conic type metrics, where $f (r) = k$ for all
$r$ large, satisfy the condtions on $f$ in the above theorem. \\[0.1in] {\bf
Proof.} \ \ If $\alpha \not\equiv 0$, if follows from the proof of lemma 2.6
 (c.f. [2]) that $\alpha > 0$ and hence 
$\alpha' (r) > 0$ for all $r > 0$. Assume that 
$$g' (r) \le a g (r)  \ \ \ \ {\mbox{for \ \ all}} \ \ r > 0\,. \leqno 
({\mbox
{I}})
$$ Then (2.22) implies that
$$
\Theta' (r) < C \Theta (r) \alpha' (r)\,,
$$  where we may take $C = C_o + 1$. Substitute the above inequality into 
(1.4)
and we have
$$ 
C\Theta [\alpha' (r)]^2 + \Theta \alpha'' (r) + (n-1)\frac{f'}{f}\Theta
\alpha'-(n-1)\frac{\Theta}{f^2} g(\alpha ) g' (\alpha ) > 0\,. \leqno (2.24)
$$
As $0 \le f' (r) \le c_o^2$, we have
$${{f' (r) }\over {f (r)}} \le c_1 \ \ \ \ {\mbox{for \ \ all}} \ \ r \ge
r_o\,,$$   where $c_1$ is a positive constant.  Since $\Theta (r) > 0$ and
$\alpha' (r) < C'
$ for all $r > r_o$, we have
$$  
{{\alpha''}\over {g(\alpha) g' (\alpha)}}  + C''\, {{\alpha'}\over {
g(\alpha) g' (\alpha)}} -(n-1)\frac{1}{f^2} > 0 \leqno (2.25)
$$ 
for all $r \ge r_o$, where $C''$ is a positive constant. For $r \ge r_o$, 
let 
$$
H (r) = {{\alpha' (r)}\over {g (y) g' (y)}}\,. \leqno (2.26)
$$ We have 
$$  {{\alpha''}\over {g(\alpha) g' (\alpha)}} = H' + H^2 [g'^2 (\alpha) 
(\alpha)
+ g (\alpha) g'' (\alpha)]\,.
$$ Substitute into (2.25) we obtain
$$  [g'^2 (\alpha) + g (\alpha) g'' (\alpha)] H^2 + H' + C'H 
-\frac{(n-1)}{f^2}
> 0.
\leqno (2.27)
$$ As $\alpha' > 0$, we have either $\lim_{r \to \infty} \alpha (r) =
\infty$ or $\alpha$ is bounded - in the latter case $\alpha \equiv 0$ by
theorem 2.17. Therefore we may assume that $\lim_{r \to \infty} \alpha 
(r) =
\infty$. Assume that 
$$g'' (r) \ge 0 \leqno ({\mbox {II}})$$ for all $r \ge r_o$. As $g' (0) = 
1$, we
have $g' (r) \ge 1$ for all $r > 0$. Thus $\lim_{r \to \infty} g (r) = 
\infty$.
Furthermore, assume that for some positive constant $c$, 
$$g'2 (r) - g (r) g'' (r) \ge c^2 \ \ \ {\mbox{for \ \ all}} \ \ r > 0\,. 
\leqno
({\mbox {III}})$$ Then we have 
\begin{eqnarray*} 0 \le [g'^2 (\alpha) + g (\alpha) g'' (\alpha)] H^2 & = &
{{\alpha'^2 [g'^2 (\alpha) + g (\alpha) g'' (\alpha)] }\over { g^2 
(\alpha) g'^2
(\alpha)}}\\ & \le & {{\alpha'^2}\over {g^2 (\alpha)}} + {{\alpha'^2 g 
(\alpha)
g'' (\alpha) }\over { g^3 (\alpha)  g''(\alpha)}}   \ \ \ \to 0 \ \ \
{\mbox{as}} \ \ \ r \
\to \
\infty\,,
\end{eqnarray*} as by assumption (III) $g'^2 (\alpha) > g (\alpha) g''
(\alpha)\,.$ Similarly we have 
$$\lim_{r \to \infty} H (r) = 0\,.
$$ 
Then (2.27) gives
$$
\lim_{r \to \infty} H' (r) > k > 0\,.
$$ 
As $H (r) > 0$, we have 
$$
\int_{r_o}^r H (s) ds = O (r)\,. \leqno (2.28)
$$ 
On the other hand by (III) we have 
$$ 
- ( \ln {{g' (y)}\over {g (y)}} )' = H [g'^2 (y) - g (y) g'' (y)] \ge c^2 H
$$  
for all $r \ge r_o$. Thus 
$$
- \ln {{g' (y (r))}\over {g (y (r))}} |^r_{r_o} = -  \int_{r_o}^r ( \ln
{{g' (y)}\over {g (y)}} )' dr \ge c^2 \int_{r_o}^r H (s) ds\,. \leqno (2.29)
$$
Assume that 
$$
{{g' (s)}\over {g (s)}} \ge a' \ \ \ {\mbox{for \ \ all}} \ \ s \ge 1\,, 
\leqno
({\mbox {IV}})
$$
where $a'$ is a positive constant, or,
$$
{{g' (s)}\over {g (s)}} \ge {{C'}\over {s}} \ \ \ {\mbox{for \ \ all}} \ 
\ s \ge
1\,. \leqno ({\mbox {IV}}')
$$
In the first case we have $a' \le g' (y)/g (y) \le a$ for all $r$ large 
enough.
Therefore (2.21) implies that 
$$
\int_{r_o}^\infty H (s) ds  < \infty\,.
$$ 
Since $\alpha' (r) < C$, which implies that $\alpha (r)
\le C_1 r$ for all
$r \ge r_o$, where $C_1$ is a positive constant. As $\lim_{r \to \infty} 
\alpha
(r) =
\infty$, we may assume that $\alpha (r) \ge 1$ for all $r$ large enough. 
Therefore in the second case we have 
$$
a \ge {{g' (y)}\over {g (y)}} \ge {{C'} \over {\alpha (r)}} \ge {{C''} 
\over r}
$$
for all $r$ large enough, where $C''$ is a positive constant. We have
$$
\int_{r_o}^r H (s) ds \le C'' \ln r
$$
for some constant $C''$. In both cases, we have a contradiction with 
(2.28).\qed
For $p = 2$, Liouville's type theorems for harmonic maps have been 
discussed in
[11]. In [12], Liouville's type theorems for $p$-harmonic maps are obtained
under assumptions on the curvature of the manifolds. The conditions in 
theorems
2.17 \& 2.23 involve the first derivative of the solution and growth of 
$f (r)$,
which relax assumptions on the curvature.

\vspace{0.5in}

{\bf \Large 3. \ \ Existence of bounded positive solutions}

\vspace{0.3in}

In this section we show the following local existence theorem for 
equation (1.4).\\[0.1in] 
{\bf Theorem 3.1. [2]} \ \ {\it Let $f$ and $g$ be $C^2$
functions on $[0,
\infty)$ which satisfy the conditions in (1.1).   Then for any $\al>0$, 
there is
a unique solution $\alpha \in C^2[0,\ep)$ to equation (1.4) for some $\ep>0$,
such that 
$\alpha (0) = 0$, $\alpha'(0)=\al$ and $\alpha > 0$ on 
$(0,\ep)$.}\\[0.1in]  {\bf Proof.} \ \
Equation (1.4) is equivalent to
$$ \alpha''+\frac{n-1}{r}\alpha'-\frac{n-1}{r^2}\alpha= (
(n-1)(\frac{1}{r}-\frac{f'}{f})-\frac{ {(\theta^{q-1})}'}{\theta^{q-1}} )\,
\alpha'+(n-1)(\frac{g(\alpha)g'(\alpha)}{f^2}-\frac{\alpha}{r})\ .
$$ Let $z=\alpha/r$. Then the above equation becomes
$$ rz''+(n+1)z'=[\frac{n-1}{r}-\frac{(n-1)f'}{f}-
\frac{ {(\theta^{q-1})}'}{\theta^{q-1}}](rz'+z)+
(n-1)(\frac{g(rz)g'(rz)}{f(r)^2}-\frac{z}{r})\,. \leqno (3.2)
$$ For equation (3.2), the homogeneous part has linearly independent 
solutions
$1$ and $r^{-n}$, with Wronskian $(-n) r^{-n-1}$. Letting
$v(r)=\sqrt{r}z'(r)$,  the second order equation (3.2) with $\lim_{r\to
0+}z(r)=\alpha$ is equivalent to the following system of integral equations:
\begin{eqnarray*} (3.3) \ \ \ \ \ z(r)&=& T_1(z,v)=\alpha+(n-1)\int_{0}^{r}
(1-\frac{s^n}{r^n})\Phi(s,v(s),z(s))\, ds \ \ \ {\mbox{and}}\\  v(r)&=&
T_2(z,v)=n(n-1)\int_0^r
\frac{s^n}{r^{n+\frac{1}{2}}}\Phi(s,v(s),z(s))\, ds\,,
\end{eqnarray*} where
$$
\Phi(s,v,z)=(\frac{1}{s}-\frac{f'(s)}{f(s)})(\sqrt{s}v+z) -\frac{
{(\theta^{q-1})}'}{(n-1)\theta^{q-1}}(\sqrt{s}v+z)
+(\frac{g(sz)g'(sz)}{f(s)^2}-\frac{z}{s})\,.
$$ 
We apply the method of successive approximation on  (3.3). Except for the
term 
$$A(s,v,z) = \frac{ {(\theta^{q-1})}'}{\theta^{q-1}}\,,$$
$\Phi(s,v,z)$ is only a polynomial in $v$ and $z$. Thus we need to make sure
that all of $\Phi$, $\frac{\partial \Phi}{\partial z}$ and
$\frac{\partial\Phi}{\partial v}$ are continuous near $s=0$ for any $z>0$ and
$v$. Consider the term $A (s, v, z)$. By (2.3), we have
\begin{eqnarray*}
A(s,v,z)&=&\frac{2(q-1)(n-1)}{(2q-1){(\sqrt{s}v+z)}^2+(n-1){g(sz)}^2
/{f(s)}^2}\, \left [
\frac{2g(sz)g'(sz)(\sqrt{s}v+z)}{{f(s)}^2} \right.\\
 & &\quad \left. -\frac{f'(s){(\sqrt{s}v+z)}^2}{ f(s)}-\frac{{g(sz)}^2
f'(s)}{{f(s)}^3}\right ]\,,
\end{eqnarray*} where the denominator is a differentiable function and never
vanishes when $s>0$. Assume that when $r$ is close to 0,
\begin{eqnarray*} f(r)&=&r +f_1r^2+o(r^2)\ ;\\ g(r)&=&r+g_1r^2+o(r^2)\ ,
\end{eqnarray*} for some $f_1,g_1\in {\real}$. A computation gives
$$
\lim_{s\to 0}A(s,v,z)=\frac{2(q-1)(n-1)}{(2q+n-2)z^2}\, (-v^2+4g_1 z^3 -4f_1
z^2)\ ,\leqno (3.4)
$$ and
\begin{eqnarray*}
\lim_{s\to 0} \frac{\partial A}{\partial v}(s,v,z)&=&
\frac{-4(q-1)(n-1)v}{(2q+n-2)z^2}=\frac{\partial A}{\partial v}(0,v,z)\ ,\\
\lim_{s\to 0} \frac{\partial A}{\partial z}(s,v,z)&=&
\frac{4(q-1)(n-1)(v^2+2g_1 z^3)}{(2q+n-2)z^3}=
\frac{\partial A}{\partial z}(0,v,z)\ .
\end{eqnarray*} It follows that $\Phi$, $\frac{\partial\Phi}{\partial z}$ and
$\frac{\partial\Phi}{\partial v}$ are all continuous near $(0,z,v)$ for any
$z>0$ and $v$. Moreover, as $s$ tends to 0,
$$
\frac{1}{s}-\frac{f'(s)}{f(s)}=f_1+O(s)\,,
\leqno (3.5)
$$ and
$$
\frac{g(sz)g'(sz)}{f(s)^2}-\frac{z}{s}=3g_1 z^2-2f_1 z+o(1)\,.
\leqno (3.6)
$$ For $\ep>0$ sufficiently small, we consider the complete metric space
$$ S=\{ (z,v):\ z,v\in C(0,\ep) \mbox{ with } \sup_{(0,\ep)}|z(r)-\al|\leq
\frac{\al}{2}
\mbox{ and } \sup_{(0,\ep)}|v(r)|\leq 1 \}\ .
$$ By equations (3.4), (3.5) and (3.6), for any
$(z,v)\in S$ and $\ep<\al^2/4$,
\begin{eqnarray*} & \ & |\Phi(s,v(s),z(s))|\\ & \leq &
|f_1|\al+\frac{16(q-1)(n-1)}{\al(2q+n-2)}
(1+\frac{27}{2}|g_1|\al^3+9|f_1|\al^2)+(\frac{27}{2}|g_1|\al+|f_1|)\al^2 
+\eta\,,
\end{eqnarray*} where $\eta>0$ is small. Thus the operator defined by
$T=(T_1,T_2)$ as given in (3.3) maps $S$ into itself. Therefore by 
contraction
mapping theorem, the system (3.3) has a unique solution $(z,v)\in 
C[0,\ep)$ for
some
$\ep>0$, with $z(0)=\al$ and $v(0)=0$. Since $z'=v/\sqrt{r}$,
$y'=\sqrt{r}v+z$, and
$$ v'(r)=\frac{(n-1)\Phi}{r}-(n+\frac{1}{2})\frac{v}{r}\,,
\leqno (3.7)
$$ we have $z,v\in C^\infty(0,\ep)$ and so $y\in C^\infty(0,\ep)\cap
C^1[0,\ep)$. It remains to show that $y$ is $C^2$ at $r=0$. To do this, we
observe that by L'Hopital's rule,
\begin{eqnarray*}
\lim_{r\to 0}z'(r) &=& \lim_{r\to 0}\frac{v}{\sqrt{r}} =\lim_{r\to 0}
(n-1)\int_0^r \frac{s^n}{r^{n+1}}\Phi(s,v(s),z(s))\, ds\\
 &=& \frac{n-1}{n+1}\Phi_0\ ,
\end{eqnarray*} which is also equal to $z'(0)$. Here
$$
\Phi_0 \stackrel{\mbox{def}}{=}\lim_{s\to 0}\Phi(s,v(s),z(s))
=\frac{(n-1)(n-2q+2)g_1\al^2+(4nq-2n^2-2q-n+2)f_1\al}{2q+n-2}\ .
$$ As a result,
$$ \alpha''(0)=\lim_{r\to 
0}\frac{\alpha'(r)-\al}{r}=\frac{2(n-1)}{n+1}\Phi_0 =\lim_{r\to
0}\alpha''(r)\ .
$$
{\bf Q.E.D.}\\[0.1in]
{\bf Theorem 3.8. [2]} \ \ {\it  For $p > 2$, let $\alpha (r) \in C^2 (0, 
R)$ be
a positive solution to (1.1) with $\lim_{r \to 0^+} \alpha (r) = 0$, where
$R$ is a positive number. Suppose that $f$ and $g$ satisfy the conditions in
(1.4) and there exists a positive constant $a$ such that $|g' (r)| \le a$ for
all 
$r \ge 0$. Then
$y$ can be extended to a positive solution of (1.4) on
$(0,
\infty)$.}\\[0.1in] {\bf Proof.} \ \ By the continuation of solutions  
[3, p.15],
it suffices to show that a solution of equation (1.4) is bounded on the 
interval
$(0, R)$. There exists a positive constant $b$  such that 
$|f' (r)| \le b$ and
$|g' (\alpha (r))| \leq a$ for all $r \in (0, R)$. By lemma 2.15 we have 
$\theta
(r) < C$ for  all
$r \in (R/2, R)$. Thus $\alpha' (r)$ and hence $\alpha (r)$ are bounded 
on $(R/2, R)$.  We
can continue the solution to $(0, R + \epsilon)$ for some positive  number
$\epsilon$. \ \ \ \ \ {\bf Q.E.D.}\\[0.2in]
{\bf Theorem 3.9. [2]} \ \ {\it Suppose that $p > 2$ and $n \ge 2$, and 
assume
that there exist constants $a>0, $ and
$\delta>1$ with
$(n-1)\delta> 2q-1$ ($q = p/2$) such that}
$$ g(y)>0,\qquad 0\leq g'(y) \leq a\quad\mbox{for any }y> 0\ ,
$$ {\it and}
$$ f(r)\asymp r^\delta\qquad f'(r)>0 \quad\mbox{for large }r\ ,
$$ {\it then all positive solutions to equation (1.1) are bounded if and 
only if
$n>1$.}\\[0.1in] {\bf Remark.} \ \ When the above conditions are 
satisfied, by
using theorem 3.1 and theorem 3.8, there exists a family of global
solutions which are all bounded harmonic maps.\\[0.1in] 
{\bf Proof.} \ \ When
$n=1$, then $\th={(\alpha')}^2$. Equation (1.4) implies that 
$\alpha'$ is constant. Hence
$$ \alpha'(r)=\alpha'(0)=\al\ ,\quad\mbox{ or }\quad \alpha(r)=\al r\ ,
$$ which is unbounded. Asume that $n>1$. Equation (1.4) implies that 
\begin{eqnarray*} {(\th^{q-1}\alpha')}'&=&
(n-1)\th^{q-1}(\frac{g(\alpha)g'(\alpha)}{f^2}-\frac{f'\alpha'}{f})\\
 &\leq& (n-1)\th^{q-1}(\frac{C\alpha}{r^{2\de}})\,,
\end{eqnarray*} when $r$ is sufficiently large, say, when $r\geq r_0$. Hence
$$ {(\th^{q-1}\alpha')}\leq C\, (\int_{r_0}^{r}s^{-2\de}\th^{q-1}\alpha
\, ds + 1)\,.
$$ 
Apply Lemma 2.15 to obtain
\begin{eqnarray*} {(\th^{q-1}\alpha)}'&\leq& 
2(n-1)(q-1)(a+1)\frac{\th^{q-1}\alpha}{r^\de}+
C\, (\int_{r_0}^r s^{-2\de}\th^{q-1}\alpha\, ds+ 1)\ ,\\
 &=&\frac{2\de}{r}\, \th^{q-1}\alpha+C(\int_{r_0}^r 
s^{-2\de}\th^{q-1}\alpha\, ds+1)
\end{eqnarray*} for $r\geq r_0$, where $r_o$ is a positive constant. Thus
$$ {(r^{-2\de}\th^{q-1}\alpha)}'
 = C r^{-2\de} \left ( \int_{r_0}^r s^{-2\de}\th^{q-1}\alpha\, ds 
+1\right )\,.\leqno
(3.10)
$$ Let
$$ u= \int_{r_0}^r s^{-2\de}\th^{q-1}\alpha\, ds + 1\,.
$$ Equation (3.10) can be written as
$$ u''\leq Cr^{-2\de}u\leq \ep^2 r^{-2} u \leqno (3.11)
$$ for some small $0<\ep<2\de-2$. So $u$ is a supersolution of the Euler
equation. By the weak maximum principle, we obtain
$$ u\leq C r^{1+\ep}\ .
$$ Substitute this result into equation (3.10) and integrate. We have
$$ r^{-2\de}\th^{q-1}\alpha\leq C\, (r^{-2\de+2+\ep}+1)\leq C
$$ for $r\geq r_0$, and so
$$
\int_{r_0}^r s^{(n-3)\de}\th^{q-1}\alpha\, ds \leq C r^{(n-1)\de+1}\ .
\leqno (3.12)
$$ On the other hand, we have
$$ {(f^{n-1}\th^{q-1}\alpha')}'\leq (n-1)af^{n-3}\th^{q-1}\alpha\,.
\leqno (3.13)
$$ Integrate from $r_0$ to $r$, we obtain
$$ f^{n-1}\th^{q-1}\alpha'\leq C\, (\int_{r_0}^r 
s^{(n-3)\de}\th^{q-1}\alpha\, ds +
1)\,.\leqno (3.14)
$$ Substitute equation (3.13) into equation (3.14) to obtain
$$ {(\alpha')}^{2q-1}\leq \th^{q-1}\alpha'\leq C r\,.
$$ After integration we have
$$ \alpha\leq C r^{1+1/(2q-1)}\,,
$$ This implies that
$$
\th^{q-1}\alpha\leq C r^2\,.
$$ Substitute this result into equation (3.14) again,
$$ {(\alpha')}^{2q-1}\leq C r^{3-2\de}+ C r^{-(n-1)\de}\ .
$$ If $2\de-2\geq (n-1)\de>2q-1$, then 
$$ \alpha'\leq C r^{\frac{-(n-1)\de}{2q-1}}\,.
$$ Therefore we have
$$ \alpha\leq C\, ( r^{1-\frac{(n-1)\de}{2q-1}}+1)\,,\leqno (3.15)
$$ which means that $\alpha$ is bounded, as required. Otherwise, 
$$ \alpha \leq C\, ( r^{1+\frac{3-2\de}{2q-1}}+ 1)\,.
$$ When $1+\frac{3-2\de}{2q-1}<0$, then $\alpha$ is bounded. If not,

$$
\th^{q-1}\alpha\leq C r^{4-2\de}\ .
$$ Substitute into equation (3.14),
$$ 
{(\alpha')}^{2q-1}\leq C r^{5-4\de}+C r^{-(n-1)\de}\ .
$$ 
In this way, we have a sequence of upper bounds for
${(y')}^{2q-1}$ as a linear combination of $r^{2k+1-2k\de}$ and 
$r^{-(n-1)\de}$.
Since $\de>0$, $r^{-(n-1)\de}$ will dominate after a finite number of
iterations. Therefore equation (3.15) holds eventually. Thus $y$ is bounded.
\ \ \ \ \ {\bf Q.E.D.}\\[0.2in]
{\bf Theorem 3.16. [2]} \ \ {\it For $p > 2$ and $n \ge 2$, if
$(n-1)\de\leq 2q-1$, where $q = p/2$,  and $f\asymp r^\de$ for large
$r$, then no positive entire solution to equation (1.4) is 
bounded.}\\[0.1in] 
{\bf Proof.} \ \  
Observe by equation (3.13), the function 
$f^{n-1}\th^{q-1}y'$ is increasing to some $P> 0$ (may be $\infty$). 
Hence there
is some constant $C$ such that for large $r$,
$$
\th^{q-1}\alpha'\geq C r^{-(n-1)\de}\ .\leqno (3.17)
$$ Suppose for contradiction that a solution $y$ is bounded. Then we must 
have
$\th\asymp {(\alpha')}^2$. For if otherwise, we have a sequence of points 
where
$$
\th\asymp r^{-2\de}\quad\mbox{and}\quad \alpha'=o(r^{-\de})\ .
$$ So
\begin{eqnarray*} f^{n-1}\th^{q-1}\alpha' &=& o(r^{\de(n-2q)})\ ,\\
  &=& o(1)\ ,
\end{eqnarray*} which is impossible, since $P>0$. Therefore by equation 
(3.17)
we have 
$$ \alpha'\geq C r^{\frac{-(n-1)\de}{2q-1}}
$$ for large $r$, so that
$y$ is greater than $C r^{1-\frac{(n-1)\de}{2q-1}}$, or $C \log r$. In both
cases, $y$ is unbounded. This gives a contradiction.{\bf Q.E.D.}\\[0.2in]
{\bf Theorem 3.18. [2]} \ \ {\it For $p > 2$ and $n \ge 2$, if
$(n-1)\de\leq 2q-1$,  $f\sim C r^\de$, and $0<g'\leq a\ (a>0)$, then any
positive entire solution $y$ of equation (1.4) satisfies}
$$ y'\sim A r^{\frac{-(n-1)\de}{2q-1}}\qquad
\mbox{for some constant } A>0\ .
$$ {\bf Proof.} \ \ As $f^{n-1}\th^{q-1}y'$ tends to a number $P$, finite or
infinite. Suppose that $(n-1)\de<2q-1$. If $P=\infty$, then we can find a
sequence $R_N\to \infty$ such that, at $r=R_N$, 
$f^{n-1}\th^{q-1}y'(R_N)=N$, and
for $r\leq R_N$,
$$ f^{n-1}\th^{q-1}y'(r)\leq N\ .
$$ 
Hence 
$ r^{(n-1)\de}{(y')}^{2q-1}\leq C N$ for some constant $C$. This gives, 
after an
integration,
$$ \alpha (r) \leq C N^{\frac{1}{2q-1}} r^{1-\frac{(n-1)\de}{2q-1}} \ .
$$ Thus
$$
\th(r)\leq C N^{\frac{2}{2q-1}}r^{-\frac{2(n-1)}{2q-1}}\ .
$$ So by equation (3.10) we have 
\begin{eqnarray*} N&=& O\, \left (\int_r^{R_N} \th^{q-1}N^{\frac{1}{2q-1}}
s^{(n-3)\de+1-\frac{(n-1)\de}{2q-1}}\, ds +1\right )\ ,\\
 &=& O \, (N\int_r^{R_N} s^{1-2\de}\, ds + 1 )\ ,\\
 &=& o(1) N
\end{eqnarray*} for $r$ sufficiently large. Therefore, as $N\to\infty$, 
$N=o(N)$ which is impossible. Hence there exists some constant $P\in {\bf R}$
such that 
$$ f^{n-1}\th^{q-1}y'\sim P\ .\leqno (3.19)
$$ Thus
$$
 \alpha'\leq C r^{\frac{-(n-1)\de}{2q-1}}\quad\mbox{and}
\quad \alpha\leq C r^{1-\frac{(n-1)\de}{2q-1}}\ .
$$ Substituting these bounds into $\th$ in equation (3.17), we obtain
$$
 y'\asymp r^{\frac{-(n-1)\de}{2q-1}}\quad\mbox{and}\quad \alpha\asymp
r^{1-\frac{(n-1)\de}{2q-1}}
$$ Therefore $\th\sim {(\alpha')}^2$ and 
$$ \alpha'\sim A r^{\frac{-(n-1)\de}{2q-1}}\quad \mbox{ for some }A>0\ .
$$
\hskip0.25in When $(n-1)\de=2q-1$, the procedure is similar.
$$ \alpha'\leq C N^{\frac{1}{2q-1}}r^{-1}\ ,
$$ and 
$$ \alpha\leq C N^{\frac{1}{2q-1}}\log r\ .
$$ So
$$ N\leq C N\int_r^{R_N} s^{1-2\de}\log s\, ds =o(1)N\ .
$$ The other cases follow in the same way.
\ \ \ \ \ {\bf Q.E.D.}\\[0.2in] 
{\bf  Corollary 3.20. [2]} \ \ {\it For $p > 2$ and
$n \ge 2$, if 
$\de>\frac{2q-1}{n-1}\geq 1$,
$f\sim C r^\de$, and $0<g'\leq a$, then} 
$$ y'\sim A r^{\frac{-(n-1)\de}{2q-1}}\quad\mbox{for some } A>0\ .
$$ {\bf Proof.} \ \ We follow the proof in Theorem 3.11. Suppose there is a
sequence
$R_N\to \infty$ such that $ f^{n-1}\th^{q-1}y'(R_N)=N$, then for
$r<R_N$,
$$ \alpha'\leq C N^{\frac{1}{2q-1}} r^{-\frac{(n-1)\de}{2q-1}}\ ,
$$ but
$$ \alpha\leq C N^{\frac{1}{2q-1}}\ .
$$ Hence 
\begin{eqnarray*} N&\leq& O\left (N\int_r^{R_N} s^{(n-3)\de}\max\{ 
s^{-(n-1)\de},
\ s^{-2\de(q-1)}\}\, ds+1\right )\, \\
 &=& O (N\int_r^{R_N} s^{-\de}\, ds +1)\ ,\\
 &=& o(1)N
\end{eqnarray*} for $r$ sufficiently large, and as $N\to\infty$.  The other
cases are similar.
\ \ \ \ \ {\bf Q.E.D.}\\[0.2in] 
{\bf  Corollary 3.21. [2]} \ \ {\it For $p > 2$ and
$n \ge 2$, if
$\de>\frac{2q-1}{n-1}\geq 1$,
$f\sim C r^\de$, and $0<g'\leq a$, then} 
$$ \alpha'\sim A r^{\frac{-(n-1)\de}{2q-1}}\quad\mbox{for some } A>0\ .
$$ {\bf Proof.} \ \ We follow the proof in Theorem 3.18. Suppose there is a
sequence $R_N\to \infty$ such that $ f^{n-1}\th^{q-1}y'(R_N)=N$, then for
$r<R_N$,
$$ \alpha'\leq C N^{\frac{1}{2q-1}} r^{-\frac{(n-1)\de}{2q-1}}\ ,
$$ but
$$ \alpha\leq C N^{\frac{1}{2q-1}}\ .
$$ Hence 
\begin{eqnarray*} N&\leq& O\left (N\int_r^{R_N} s^{(n-3)\de}\max\{ 
s^{-(n-1)\de},
\ s^{-2\de(q-1)}\}\, ds+1\right )\, \\
 &=& O (N\int_r^{R_N} s^{-\de}\, ds +1)\ ,\\
 &=& o(1)N
\end{eqnarray*} for $r$ sufficiently large, and as $N\to\infty$.  The other
cases are similar.\ \ \ \ \ {\bf Q.E.D.}\\[0.2in]
{\bf Theorem 3.22. [9]} \ \ {\it Assume that $f' (r) > 0$ for
all
$r > 0$ and $\lim_{r
\to
\infty} f' (r) = \infty$. Suppose that
$g' (y)
\le k$ for all $y \ge 0$, where $k$ is a positive constant. Let
$\alpha (r)
\in C^2 (0,
\infty)$ be a positive solution to (1.4) with $f$ and $g$ as above and 
$\lim_{r
\to 0^+}
\alpha (r) = 0$. If $2 < p \le n$, then for any $\epsilon \in (0, 1)$, there
exist positive constants $r_o$ and $C$ such that $\alpha' (r) \le C/ f^{1 -
\epsilon} (r)$ for all $r \ge r_o$. If $p > n \ge 2$, then for any 
$\epsilon >
0$, there exist positive constants $r_o$ and $C$ such that} 
$$\alpha' (r) \le { C \over { (f (r))^{(1 - \epsilon) {{(n - 1)}\over {(p 
-1)}}}
}}$$  {\it for all $r \ge r_o$.}\\[0.2in]   {\bf Proof.} \ \ Given $1 > 
\epsilon
> 0$, as in (5.15) (with $\epsilon = 1 -
\tau$) we can find a positive number
$r_o$ such that
\begin{eqnarray*} & \ & \ {{2 g (\alpha (r)) g' (\alpha (r))
\alpha' (r)}\over {f^2 (r)}} - {{f' (r)}\over {f (r)}} [ {{ g^2 (\alpha
(r))}\over {f^2 (r)}} + (\alpha' (r))^2 ]\\  & \le & {{ g^2 (\alpha 
(r))}\over
{\epsilon f^3 (r) f' (r)}} [k^2 - \epsilon^2 f' (r)^2] - (1 - \epsilon) 
[{{ f'
(r) g^2 (\alpha (r))}\over {f^3 (r)}} + {{f' (r)}\over {f (r)}} (\alpha'
(r))^2]\\ &
\le & - (1 - \epsilon) {{f' (r)}\over {f (r)}} [{{ \alpha^2 (r)}\over 
{f^2 (r)}}
+ (\alpha' (r))^2]
\end{eqnarray*}  for all $r \ge r_o$, as $g' (s) \le k$ and $\lim_{r \to 
\infty}
f' (r) =
\infty$. By (2.15) we have
\begin{eqnarray*} (3.23) \ \ \ \ \ \ \ \ \ & \ & \ 
\Theta' (r) [ ( p - 1 ) (\alpha' (r))^2 + (n - 1) {{g^2 (\alpha (r) 
)}\over {f^2
(r)}} ]\\ & \le & - (1 - \epsilon) (p - 2) (n - 1) {{f' (r)}\over {f 
(r)}} [{{
g^2 (\alpha (r))}\over {f^2 (r)}} + (\alpha' (r))^2] \Theta (r) \ \ \ \ \ \
\ \ \ \ \ \ \ \ \ \ \ \ \ \ \  \ \ \ \ \ \ \ \ \ \ \ 
\end{eqnarray*} for all
$r \ge r_o$.  Thus if $2 < p \le n$, then (3.23) implies that  
$$\Theta' (r) \le - (1 - \epsilon) (p - 2) {{f' (r)}\over {f (r)}} \Theta 
(r)$$
for all $r \ge r_o$. An integration gives $\alpha' (r) \le C/ f^{1 - 
\epsilon}
(r)$ for some positive constant $C$ and for all $r \ge r_o$. If $p > n > 
2$, then
(3.19) gives 
$$\Theta' (r) \le - (1 - \epsilon) {{(p - 2)(n - 1)}\over {p - 1}} {{f' 
(r)}\over
{f (r)}}
\Theta (r)$$ for all $r \ge r_o$.  An integration gives
$$\alpha' (r) \le {C \over {(f (r))^{(1 - \epsilon) {{(n - 1)}\over {(p 
-1)}}}
}}$$ for some positive constant $C$ and for all $r \ge r_o$. \ \ \ \ \ {\bf
Q.E.D.}\\[0.2in]  
\hspace*{0.5in}As a corollary, we have the following result (c.f. 
[2]).\\[0.1in]
{\bf Corollary 3.24. [9]} \ \ {\it Assume that $f' (r) > 0$  for all $r > 
0$ and
$\lim_{r \to
\infty} f' (r) = \infty$. Suppose that
$g' (y) \le k$ for all $y > 0$, where $k$ is a positive constant. Assume that
$f (r)
\ge C r^s$ for some positive constant $C$ and for all $r > \bar r > 0$, where
$s > 1$ if $2 < p \le n$, and} 
$$s > {{(p - 1)}\over {(n - 1)}}$$  {\it if $p > n \ge 2$. Let
$\alpha (r)
\in C^2 (0,
\infty)$ be a positive solution to (1.1) with $f$ and $g$ as above and 
$\lim_{r
\to 0^+}
\alpha (r) = 0$. Then $\alpha$ is a bounded function on 
${\real}^+$.}\\[0.2in]
{\bf Proof.} \ \ From theorem 3.21, we can choose $r_o > 0$ and $\epsilon 
> 0$ so
that 
$$\alpha' (r) \le {C\over {r^{1 + \delta} (r)}}$$ for all $r \ge r_o$. Here
$\delta$ is a positive constant. An integration shows that $\alpha$ is a 
bounded
function on ${\real}^+$.\qed In particular, for $p > 2$, all rotationally
symmetric $p$-harmonic maps from the hyperbolic space to the Euclidean 
space are
bounded.

\pagebreak

{\bf \Large 4. \ \ $p$-harmonic maps from hyperbolic space to itself}

\vspace{0.4in}

In this section we discuss rotationally symmetric $p$-harmonic maps from the
hyperbolic space to itself. In this case $f (r) = g (r) = \sinh r$ as the 
metric
$$dr^2 + \sinh^2 r \, d\varphi^2$$ is exactly the hyperbolic metric.
The results in this section can be found in [8].\\[0.2in] 
{\bf Lemma 4.1.} \ \
{\it For
$p > 2$, let
$F (r,
\varphi) = (\alpha (r),
\varphi)$ be a rotationally symmetric $p$-harmonic map from the 
hyperbolic space
to itself, where $\alpha \in C^2 (0, \infty)$ and $\lim_{r \to 0^+}
\alpha (r) = \alpha (0) = 0$. If there exists $r_o > (\ln 3) /2$ such that
$\alpha (r_o) = r_o + \delta$ for some
$\delta \ge 0$ and
$\alpha' (r_o) > 1$, then
$\alpha (r) > r + \delta$ and $\alpha' (r) > 1$ for all $r > r_o$. If there
exists
$r_o > (\ln 3) /2$ such that
$\alpha (r_o) = r_o - \delta$ for some $\delta \ge 0$ and $\alpha' (r_o) 
< 1$,
then $\alpha (r) <  r_o - \delta$ and $\alpha' (r) < 1$ for all
$r > r_o$.}\\[0.2in] {\bf Proof.} \ \ Suppose that there exists $r_o > 
(\ln 3)
/2$ such that
$\alpha (r_o) = r_o + \delta$ for some
$\delta \ge 0$ and $\alpha' (r_o) > 1$. Assume that there is a point $r' 
> r_o$
such that $\alpha (r') = r' + \delta$. Then there exists a point $r \in (r_o,
 r')$ such that $\alpha (r) > r + \delta$, $\alpha' (r) = 1$ and 
$\alpha'' (r)
\le 0$. Let $\alpha (r) = r + \epsilon$ where $\epsilon > \delta$ is a 
constant.
As $f (r) = g (r) = \sinh r$, at the point $r$ where $\alpha' (r) = 1$ we 
have
\begin{eqnarray*} & \ & \ \ {{2 g (\alpha (r)) g' (\alpha (r)) \alpha' 
(r)}\over
{f^2 (r)}}  - {{f' (r)}\over {f (r)}} [ {{g^2 (\alpha (r))}\over {f^2 
(r)}} +
(\alpha' (r))^2 ]\\ & = & {1\over 8 f^3 (r)} [ 2 (e^r - e^{-r}) (e^\alpha -
e^{-\alpha}) (e^\alpha + e^{-\alpha}) \\ & \ & \ \ \ \ \ \ \ \ \ \ \ \ \  
- (e^r
+ e^{-r}) (e^\alpha - e^{-\alpha})^2 - (e^r + e^{-r}) (e^r - e^{-r})^2 
]\\  & =
& {1\over 8 f^3 (r)} (e^r e^{2\alpha} - 3e^r e^{- 2\alpha} - 3 e^{-r}
e^{2\alpha} + e^{-r} e^{-2\alpha} + 3 e^r + 3 e^{-r} - e^{-3r} - e^{3r})\,.
\end{eqnarray*}  As $\alpha (r) = r + \epsilon$, we have
\begin{eqnarray*} (4.2) \ \ \ \ \ \ \ \ \ \ \ \  & \ & \ \ {{2 g (\alpha 
(r)) g'
(\alpha (r))
\alpha' (r)}\over {f^2 (r)}}  - {{f' (r)}\over {f (r)}} [ {{g^2 (\alpha
(r))}\over {f^2 (r)}} + (\alpha' (r))^2 ]\\ & = & {1\over 8 f^3} [
(e^{2\epsilon} - 1) ( e^{3r} - 3 e^r) + (3 e^{-r} - e^{-3r}) (1 -
e^{-2\epsilon}) ]\,. \ \ \ \ \ \ \  \ \ \ \ \ \ \ 
\end{eqnarray*} If $r > (\ln 3) /2$, then $e^{3r} - 3 e^r > 0$ and $3 
e^{-r} -
e^{-3r} > 0$. Equation (4.2) implies that 
$${{2 g (\alpha (r)) g' (\alpha (r)) \alpha' (r)}\over {f^2 (r)}}  - {{f'
(r)}\over {f (r)}} [ {{g^2 (\alpha (r))}\over {f^2 (r)}} + (\alpha' 
(r))^2 ] >
0$$  and by (2.3) $\Theta' (r) > 0$ because $\theta (r) > 0$. As 
$$\Theta' (r) = ({p\over 2} - 1) \theta^{{p\over 2} - 2} (r) \theta' (r)\,,$$
hence $\theta' (r) > 0$ if $p > 2$. We have
$$\theta' (r) = (n - 1) ( {{g^2 (\alpha)} \over {f^2 (r)}} )' + 2 \alpha' (r)
\alpha'' (r)\,.$$ And
\begin{eqnarray*}
(4.3) \ \ \ \ \ \ \ \ \ \ \ \ \ \ \ \ \ \ \ \ \ \ ( {{g^2 (\alpha)} \over 
{f^2
(r)}} )' & = & {{2 f (r) g (\alpha) [ f(r) g' (\alpha) \alpha' (r) - f' 
(r) g
(\alpha) ]}\over {f^4 (r)}}\\ & = & {{g (\alpha)}\over {f^3 (r)}} (e^r
e^{-\alpha} - e^{-r} e^{\alpha} )\,. \ \ \ \ \ \ \ \ \ \ \ \ \ \ \ \ \ \ 
\ \ \ \
\ \ \ \ \ \ 
\end{eqnarray*} As $\alpha (r) > r$, we have 
$$( {{g^2 (\alpha)} \over {f^2 (r)}} )' < 0\,.$$  Hence $\alpha'' (r) > 
0$. This
is a contradiction. Hence $\alpha (r) > r + \delta$ for all $r > r_o$. 
Suppose
that there is a point $r'$ such that
$\alpha (r') = 1$, then there is a point $r \in (r_o, r']$ such that 
$\alpha (r)
> r + \delta$, $\alpha' (r) = 1$ and $\alpha'' (r) \le 0$. This would 
give us
the same contradiction. Hence $\alpha' (r) > 1$ for all $r > r_o$. Similar
argument works for the case $\alpha (r_o) = r_o - \delta$. \ \
\ \ \ {\bf Q.E.D.}\\[0.2in] {\bf Lemma 4.4.} \ \ {\it For $p > 2$, let $F (r,
\varphi) = (\alpha (r),
\varphi)$ be a rotationally symmetric $p$-harmonic map from the 
hyperbolic space
to itself, where $\alpha \in C^2 (0, \infty)$ and $\lim_{r \to 0^+}
\alpha (r) = \alpha (0) = 0$. Then we can find a point $r_o > 0$ such that
either $\alpha (r) \ge r$ or $\alpha (r) \le r$ for all
$r > r_o$.}\\[0.2in] {\bf Proof.} \ \ Suppose that we cannot find such a 
point
$r_o$. Then there exist two points $r_2 > r_1 > (\ln 3)/2$ such that $\alpha
(r_1) = r_1$ and
$\alpha (r_2) = r_2$ and $\alpha (r) < r$ or $\alpha (r) > r$ for all $r \in
(r_1\,, r_2)$. Then we can find a point $r \in (r_1\,, r_2)$ such that 
$\alpha' 
(r) = 1$ and if $\alpha (r) > r$, then $\alpha'' (r) \le 0$; if $\alpha 
(r) <
r$, then $\alpha'' (r) \ge 0$. But the proof of lemma 4.1 shows that 
these are
impossible. \ \ \ \ \ {\bf Q.E.D.}\\[0.2in]
\hspace*{0.5in}We have the following Liouville's type theorem for 
$p$-harmonic
map on hyperbolic spaces.\\[0.2in] {\bf Theorem 4.5.} \ \ {\it For $p > 
2$, let
$F (r, \varphi) = (\alpha (r),
\varphi)$ be a rotationally symmetric $p$-harmonic map from the 
hyperbolic space
to itself, where $\alpha \in C^2 (0, \infty)$ and $\lim_{r \to 0^+}
\alpha (r) = \alpha (0) = 0$. Assume that  there exists $r_o > (\ln 3) 
/2$ such
that $\alpha (r_o) > r_o$ and $\alpha' (r_o) > 1$. Then given any positive
number $C$, we have a point $r_c$ such that $\alpha (r) > r + C$ for all 
$r >
r_c$.}\\[0.2in]  {\bf Proof.} \ \ Assume that there is a positive number $C$
such that $\alpha (r)
\le r + C$ for all $r > r_o$. Let $\alpha (r_o) = r_o + \epsilon_o$ for some
constant
$\epsilon_o > 0$. By lemma 4.1, we have $\alpha (r) > r + \epsilon_o$ and
$\alpha' (r) > 1$ for all
$r \ge r_o$. Assume that at a point $r >\bar r$, $\alpha (r) = r + 
\epsilon$ and
$\alpha' (r)
\le 1 + \eta$, where $\epsilon \ge \epsilon_o$, and
$\bar r > r_o$ and $\eta$ are a positive constants to be determined 
later. Then
\begin{eqnarray*} & \ & {{2 g (\alpha (r)) g' (\alpha (r)) \alpha' (r)}\over
{f^2 (r)}}  - {{f' (r)}\over {f (r)}} [ {{g^2 (\alpha (r))}\over {f^2 
(r)}} +
(\alpha' (r))^2 ]\\ & \ge & {1\over {8 f^3 (r)}} [ 2 (e^{3r} e^{2 
\epsilon} -
e^{-r} e^{-2
\epsilon} - e^r e^{2 \epsilon} + e^{-3r} e^{-2 \epsilon})\\ & \ & \ \ \ \ 
- 
(e^{3r} e^{2 \epsilon} - 2 e^r + e^{-r} e^{-2 \epsilon} + e^{r} 
e^{2\epsilon}  -
2 e^{-r} + e^{-3r} e^{-2 \epsilon})\\ & \ & \ \ \ \ - (1 + \eta)^2 
(e^{3r} - e^r
- e^{-r} + e^{-3r})]\\ & = & {1\over { (e^r - e^{-r})^3}} \,\{ e^{3r} [e^{2
\epsilon} - (1 +
\eta)^2]\\ & \ & \ \ \ \ - e^r [3e^{2\epsilon} - 2 + (1 +
\eta)^2] - e^{-r}  [3e^{2\epsilon} - 2 + (1 +
\eta)^2]\\ & \ & \ \ \ \  - e^{-3r} [e^{-2\epsilon} - (1 + \eta)^2] \}\,. 
\end{eqnarray*}  Therefore if we choose $\bar r > r_o$ to be large and fix
$\eta$ such that 
$$e^{2 \epsilon_o} > (1 + \eta)^2 + {{(e^{2 \epsilon_o} - 1)}\over 2}\,,
\leqno (4.6)$$  then there exists a positive constant
$c > 0$ such that 
$${{2 g (\alpha) g' (\alpha)
\alpha'}\over {f^2}} - {{f'}\over {f}} [ {{ g^2 (\alpha)}\over {f^2 (r)}} +
(\alpha')^2 ] \ge c$$ whenever $r > \bar r$ and $\alpha' (r) \le 1 + 
\eta$. That
is,
$$\Theta' (r) [ (p - 1) (\alpha' (r))^2 + (n - 1) {{g^2 (\alpha)}\over {f^2
(r)}}] \ge c \Theta (r)$$ for all $r > \bar r$ with $\alpha' (r) \le 1 + 
\eta$.
As $\alpha' (r) \le 1 +
\eta$ and $\alpha (r) \le r + C$ for some positive constant $C$, we have 
\begin{eqnarray*}
 (4.7) \ \ \ \ \ \ \ \ \ \ \ & \ & \ (p - 1) (\alpha' (r))^2 + (n - 1) {{g^2
(\alpha)}\over {f^2 (r)}}\\ & \le & (p - 1) (1 + \eta)^2 + {{(n - 
1)}\over 4} {
{e^{2 C} - 2 e^{-2r} + e^{-4r} e^{-2C} }\over { 1 - 2 e^{-2r} + e^{-4r} 
}} \le
c'\,, \ \ \ \ \ \ \ \ 
\end{eqnarray*} where $c'$ is a positive constant. Hence 
$$\Theta' (r) \ge c_o \Theta (r) \ \ \ \ {\mbox {for \ all}} \ \ \ r > 
\bar r \ \
\ {\mbox {with}} \ \ \ \ \alpha' (r) \le 1 + \eta\,, \leqno (4.8)$$  
where $c_o =
c/c'$ is a positive constant. Since the set $O = \{ r > \bar r  | \
\alpha' (r) > 1 + \eta \}$ is an open set, we can decompose $O$ as
$$O = \cup_{i = 1 }^{\beta} (a_i, b_i)\,,$$ where $\beta \in {\bf N} \cup \{
\infty \}$ and $a_i < b_i < a_{i + 1}$. For $r
\in(a_i, b_i)$, we have $\alpha' (r) > 1 + \eta$ and
$$( {{g^2 (\alpha)} \over {f^2 (r)}} )' \ge {{f (r) g^2 (\alpha)}\over {2 f^4
(r)}} [ \eta e^r e^\alpha + (2 + \eta) e^r e^{-\alpha} - (2 + \eta) e^{-r}
e^\alpha - \eta e^{-r} e^{-\alpha} ] \ge 0$$ if we choose $\bar r$ to be 
large.
Since $\alpha' (a_i) = \alpha' (b_i) = 1 +
\eta$ and on $(a_i, b_i)$ we have
$$( {{g^2 (\alpha)} \over {f^2 (r)}} )' \ge 0\,,$$   therefore 
$$\theta (a_i) = (n - 1) {{g^2 (\alpha (a_i))}\over {f^2 (a_i)}} + [\alpha'
(a_i)]^2 \le (n - 1) {{g^2 (\alpha (b_i))}\over {f^2 (b_i)}}   + [\alpha'
(b_i)]^2 = \theta (b_i)\,.$$ That is,
$$\Theta (a_i) \le \Theta (b_i) \leqno (4.9)$$ for all $i = 1, 2,..., 
\beta$. On
the open set $O$ we have $\alpha' > 1 +
\eta$. As $\alpha (r) > 1$ for all $r > r_o$, by lemma 4.1, we have  
$$r + C \ge \alpha (r) \ge r + \eta [\sum_{i = 1}^{\beta (r)} (b_i - 
a_i)]\,,$$
where $\beta (r)$ is the biggest index $j$ such that $b_j < r$. Therefore 
$$C \ge \eta \sum_i^{\beta} (b_i - a_i)\,.$$ So the compliment of $O$ on 
$(a_1,
\infty)$, which is given by $U_{i = 1}^{\beta} [b_i, a_{i + 1}]$,  has 
infinite
measure, where
$a_{\beta + 1} = \infty$ if $\beta$ is an finite integer. For
$r
\in [b_i, a_{i + 1}]$ we have $\alpha' (r) \le 1 + \eta$ and by (4.5), 
$\Theta'
(r)
\ge c_o
\Theta (r)$, or 
$$\Theta (a_{i + 1}) \ge \Theta (b_i) e^{c_o (a_{i + 1} - b_i)}$$ for all 
$i$.
Apply (4.6) we obtain
\begin{eqnarray*}
\Theta (a_{i + 1}) & \ge & \Theta (b_i) e^{c_o (a_{i + 1} - b_i)}\\ & \ge &
\Theta (a_i) e^{c_o (a_{i + 1} - b_i)}\\ & \ge & \Theta ( b_{i - 1}) 
e^{c_o (a_i
- b_{i-1})}e^{c_o (a_{i + 1} - b_i)}\\ & = & \Theta ( b_{i - 1}) e^{ c_o 
(a_{i +
1} - b_i) + (a_i - b_{i-1})}\,.
\end{eqnarray*} Therefore by induction we have
$$\Theta (a_{N + 1}) \ge \Theta (b_1) e^{c_o \sum_{i = 1}^N (a_{i + 1} -
b_i)}\,. \leqno (4.10)$$ For $r = a_{N + 1}$ we have $\alpha' (a_N) \le 1 
+ \eta$
and
$\alpha (a_N)
\le a_N + C$, therefore 
$$\Theta (a_N) \le c_1$$ as in the proof of lemma 4.4, where $c_1$ is a 
positive
constant independent of
$N$. This is a contradiction, as we have that the set 
$U_{i = 1}^{\beta} [b_i, a_{i + 1}]$ has infinite measure. Therefore 
there is a
point $r_c$ such that $\alpha (r_c) > r_c + C$. Lemma $4.1$ shows that 
$\alpha
(r) > r + C$ for all $r > r_c$. \ \ \ \ \ {\bf Q.E.D.}\\[0.2in]
{\bf Theorem 4.11.} \ \ {\it For $p > 2$, let $F (r, \varphi) = (\alpha (r),
\varphi)$ be a rotationally symmetric $p$-harmonic map from the 
hyperbolic space
to itself, where $\alpha \in C^2 (0, \infty)$ and $\lim_{r \to 0^+}
\alpha (r) = \alpha (0) = 0$. Assume that  there exists $r_o > (\ln 3) 
/2$ such
that $\alpha (r_o) < r_o$ and $\alpha' (r_o) < 1$. Then $\alpha$ is 
bounded on
$[0, \infty)$.}\\[0.2in] {\bf Proof.} \ \ It follows from lemma 4.1 that 
there
is a positive constant
$\epsilon_o > 0$ such that $\alpha (r) \le r - \epsilon_o$ for all $r > 
r_o$. At
a point $r > r_o$, let $\alpha (r) = r - \epsilon$ where $\epsilon \ge
\epsilon_o$. We have 
\begin{eqnarray*} & \ & \ {{2 g (\alpha (r)) g' (\alpha (r))
\alpha' (r)}\over {f^2 (r)}} - {{f' (r)}\over {f (r)}} [ {{ g^2 (\alpha
(r))}\over {f^2 (r)}} + (\alpha' (r))^2 ]\\ & = & {{2 g (\alpha (r)) g' 
(\alpha
(r))}\over {f^{3\over 2} (r)f'^{1\over 2} (r)}} [ {{f'(r) (\alpha' 
(r))^2}\over
{f (r)}} ]^{1\over 2} - {{ f' (r) g^2 (\alpha (r))}\over {f^3 (r)}} - {{f'
(r)}\over {f (r)}} (\alpha' (r))^2\\ & \le & {1\over {1 - \tau}} {{g^2 
(\alpha
(r)) g'^2 (\alpha (r))}\over {f^3 (r) f' (r)}} + (1 - \tau) {{f' 
(r)}\over {f
(r)}} (\alpha' (r))^2 - {{ f' (r) g^2 (\alpha (r))}\over {f^3 (r)}}\\ & \ 
& \ \
\ \ \ \ \ \ - {{f' (r)}\over {f (r)}} (\alpha' (r))^2\\ & = & {1\over {1 -
\tau}} {{g^2 (\alpha (r)) g'^2 (\alpha (r))}\over {f^3 (r) f' (r)}} - ( 
1- \tau)
{{ f' (r) g^2 (\alpha (r))}\over {f^3 (r)}}\\ & \ & \ \ \ \ \ \ \ \  - 
\tau [{{
f' (r) g^2 (\alpha (r))}\over {f^3 (r)}} + {{f' (r)}\over {f (r)}} (\alpha'
(r))^2]\\  & = & {{g^2 (\alpha (r))}\over {4 f^3 (r) f'(r)}} [ 
{{e^{2\alpha (r)}
+ 2 + e^{-2 \alpha (r)}}\over {1 - \tau}} - ( 1- \tau)(e^{2r} + 2 + 
e^{-2r}) ]\\ 
& \ & \ \ \ \ \ \ \ \ - \tau [{{ f' (r) g^2 (\alpha (r))}\over {f^3 (r)}} 
+ {{f'
(r)}\over {f (r)}} (\alpha' (r))^2]\\ & = & {{g^2 (\alpha (r))}\over {4 
f^3 (r)
f' (r) }} ( {{e^{2r} e^{-2\epsilon} + 2 + e^{-2r} e^{2
\epsilon} -(1 - \tau)^2 (e^{2r} + 2 + e^{-2r})}\over {1 - \tau}} )\\  & \ 
&  \ \
\ \ \ \ \ \ - \tau [{{ f' (r) g^2 (\alpha (r))}\over {f^3 (r)}} + {{f' 
(r)}\over
{f (r)}} (\alpha' (r))^2]\\ & = & - {{g^2 (\alpha (r)) e^{2r}}\over {4 
f^3 (r)
f' (r)}} [ {{ (1 - \tau)^2 - e^{-2
\epsilon} - 2e^{-2r} [(1 - \tau)^2 + 1] - e^{-4r} (e^{2\epsilon} - (1 -
\tau)^2)}\over {1 -
\tau}}]\\ & \ &  \ \ \ \ \ \ \ \ \ - 
\tau [{{f' (r) g^2 (\alpha (r))}\over {f^3 (r)}} + {{f' (r)}\over {f (r)}}
(\alpha' (r))^2]\\  & \le & - {{g^2 (\alpha (r)) e^{2r}}\over {4 f^3 (r) f'
(r)}} [ {{ (1 - \tau)^2 - e^{-2
\epsilon_o} - 2e^{-2r} [(1 - \tau)^2 + 1] -  e^{-2r} + e^{-4r} (1 -
\tau)^2 }\over {1 - \tau}}]\\  & \ & \ \ \ \ \ \ \ \ \ \ -
\tau [{{f' (r) g^2 (\alpha (r))}\over {f^3 (r)}} + {{f' (r)}\over {f (r)}}
(\alpha' (r))^2]\,,\\ 
\end{eqnarray*} where we have made use of the inequality $2AB \le A^2/(1 
- \tau)
+ (1 -
\tau)B^2$ for $1 > \tau > 0$, and $r \ge
\epsilon \ge \epsilon_o > 0$. Choose $\tau > 0$ such that $(1 - \tau)^2 
\ge e^{-2
\epsilon_o} + (1 - e^{-2\epsilon_o})/2 = (1 + e^{-2 \epsilon_o})/2$. 
Since $g'
(r) = f' (r) = (e^r + e^{-r} )/2$, we can find positive constants $\bar 
r$, such
that for $r > \bar r$, we have
$${{g^2 (\alpha (r)) e^{2r}}\over {4 f^3 (r) f' (r)}} [ {{ (1 - \tau)^2 - 
e^{-2
\epsilon_o} - 2e^{-2r} [(1 - \tau)^2 + 1] -  e^{-2r} + (1 - \tau)^2)
e^{-4r}}\over {1 - \tau}}] > 0\,. \leqno (4.12)$$ Hence for all $r > \bar r$,
\begin{eqnarray*} & \ & {{2 g (\alpha (r)) g' (\alpha (r))
\alpha' (r)}\over {f^2 (r)}} - {{f' (r)}\over {f (r)}} [ {{ g^2 (\alpha
(r))}\over {f^2 (r)}} + (\alpha' (r))^2 ]\\ & \le & - 
\tau [{{f' (r) g^2 (\alpha (r))}\over {f^3 (r)}} + {{f' (r)}\over {f (r)}}
(\alpha' (r))^2]\\ & = & - \tau {{f' (r)}\over {f (r)}} [ {{g^2 (\alpha
(r))}\over {f^2 (r)}} + (\alpha' (r))^2]\\ & \le & - c [ {{g^2 (\alpha
(r))}\over {f^2 (r)}} + (\alpha' (r))^2]\,,
\end{eqnarray*} where $c > 0$ is a positive constant, as $f' (r) /f (r)$ 
is a
bounded function for $r > \bar r > 0$. Hence
$$\Theta' (r) [ (p - 2) (\alpha' (r))^2 + (n - 1) {{g^2 (\alpha)}\over {f^2
(r)}}] \le - c\, [ {{g^2 (\alpha (r))}\over {f^2 (r)}} + (\alpha' (r))^2] 
\Theta
(r)\,. \leqno (4.13)$$  If we take $b = \max \ \{ p - 2, n - 1 \}$ and 
$c' =
c/b$, then
$$\Theta' (r) \le - c' \Theta (r) \ \ \ \ {\mbox {for \ all}} \ \ \ r > \bar
r\,. \leqno (4.14)$$ Hence 
$$\Theta (r) \le e^{-c' (t - \bar r)} \Theta (\bar r)\,.$$ So we have
$$\alpha' (r) \le C' e^{- {{c' (t - \bar r)}\over {p - 2}}} \ \ \ \ 
{\mbox {for \
all}} \ \ \ \ r > \bar r\,,$$  where $C' = (\Theta (\bar r))^{1\over {p - 
2}}$
is a positive constant unless
$\alpha$ is the trivial solution. As $\alpha' (r) \ge 0$, upon 
integration we
have that $\alpha$ is bounded. \ \ \ \ \ {\bf Q.E.D.}\\[0.2in] {\bf 
Remark.} \ \
As $\bar r$ can be estimated by using $\epsilon_o$ and $r_o$ and thus $\Theta
(\bar r)$ can be estimated in terms of $\epsilon_o$ and $r_o$. Therefore 
we can
find a constant $C (r_o, \epsilon_o)$ such that $\alpha (r)
\le C (r_o, \epsilon_o)$ for all $r > 0$. As we have $\alpha' (r) > 0$ 
for all $r
> 0\,,$ if $\alpha$ is not a trivial solution, hence under the assumption of
theorem 4.1, there is a positive number $a$ such that $\lim_{r \to \infty}
\alpha (r) = a$.\\[0.2in]  
{\bf Theorem 4.15.} \ \ {\it For $p > 2$, let $F (r,
\varphi) = (\alpha (r),
\varphi)$ be a rotationally symmetric $p$-harmonic map from the 
hyperbolic space
to itself, where $\alpha \in C^2 (0, \infty)$ and $\lim_{r \to 0^+}
\alpha (r) = \alpha (0) = 0$. Suppose that there exist positive constants $C$
and $\bar r$ such that $r + C \ge \alpha (r) \ge r - C$ for all $r > \bar r$,
then $\alpha$ is asymptotic to the the identity map $\alpha_I (r) = 
r$}\\[0.2in] 
{\bf Proof.} \ \ By lemma 4.4, we can find a constant
$r_1$ such that either $\alpha (r) \ge r$ or $\alpha (r) \le r$ for all 
$r >
r_1$. Assume that there is a point $r_o > {\mbox {max}} \{ \bar r\,, r_1 \}$
such that
$\alpha (r_o) > r_o$. If there is a point $r_2 > r_o$ such that $\alpha 
(r_2) =
r_2$, then either $\alpha (r) = r$ for all $r \ge r_2$ or there is a 
point $r_3$
such that $\alpha (r_3) > r_3$. In the second case then we can find a point
$r_4$ such that $r_2 < r_4 \le r_3$ and $\alpha (r_4) > r_4$ and
$\alpha' (r_4) > 1$. By theorem 4.5, $\alpha (r) > r + C$ when $r$ is large.
Therefore $\alpha (r) > r$ for all $r > r_o$ if $\alpha (r) \not\equiv r$ on
$(r_2, \infty)$. Let 
$$c_o = \inf \,\,\{ c > 0\  | \ {\mbox {there \ exists \ a \ point}} \ 
r_c > r_o
\
\ {\mbox {such \ that \ }} \ \alpha (r_c) = r_c + c \}\,.$$ If $\alpha 
(r_c) =
r_c + c$ for some $r_c > r_o$, then $\alpha (r) \le r + c$ for all $r > r_c$.
For if there is a point $r_5$ such that $r_5 > r_o$ and
$\alpha (r_5) > r_5 + c$, then we can find $r'_5 \in (r_c\,, c_5]$ such that
$\alpha (r'_5) > r'_5$ and $\alpha' (r_5) > 1$. By theorem 4.5 this is not
possible in this case. Given any $\epsilon > 0$, we have a point $r_{c_o +
\epsilon}$ such that $\alpha (r_c) = r_{c_o + \epsilon} + c_o + \epsilon$ and
$\alpha (r) \le r + c_o + \epsilon$ for all $r \ge r_{c_o + \epsilon}$. As
$\alpha (r) > r + c_o$ for all $r > r_o$, we have the line $y = x + c_o$ 
being
an asymptotic of $\alpha$. If $c_o > 0$, then we have $\alpha (r) > r + 
c_o$ for
all $r > r_o$. Theorem 4.5 shows that $\alpha' (r) \le 1$ for all $r > 
r_o$. In
fact if there is a point $r > r_o$ such that $\alpha' (r) = 1$, then as 
in (4.2)
and the proof of lemma 4.1, we have
$\alpha'' (r) > 0$. But this is not possible by theorem 4.5. Thus 
$\alpha' (r) < 1$ for all $r > r_o$. Since $\alpha (r) \ge r + c_o$ for 
all $r >
r_o$, we can find sequence $r_i > r_o$, $i \in {\bf N}$ with $\lim_{i \to
\infty} r_i = \infty$, such that
$\alpha (r_i) > r_i + c_o$ with $c_o > 0$ and $1 - \eta \le \alpha' (r_i) 
< 1$,
where we can choose
$\eta$ as small as we like. We choose $\eta > 0$ so small such that
$2e^{2c_o} (1 - \eta) > 2 + (e^{2\epsilon_o} - 1)/2$, that is, $3
e^{2\epsilon_o} - 4 \eta > 3$. Then for any point
$r$ with
$\alpha' (r)
\ge 1 -
\eta$ and
$\alpha (r) = r + c'_o$ with $C \ge c'_o \ge c_o$, a calculation as in (4.2)
shows that we can find a constant
$r_6 > r_o$ such that if $r > r_6$, 
\begin{eqnarray*}(4.16) \ \ \ \ \ \ \ \ & \ & \ {{2 g (\alpha (r)) g' (\alpha
(r))
\alpha' (r)}\over {f^2 (r)}} - {{f' (r)}\over {f (r)}} [ {{ g^2 (\alpha
(r))}\over {f^2 (r)}} + (\alpha' (r))^2 ]\\  & \ge & {1\over {8 f^3 (r)}} 
[ 2 (
1 - \eta) (e^{3r} e^{2 c_o} - e^{-r} e^{-2 c_o} - e^r e^{2 c'_o} + 
e^{-3r} e^{-2
c'_o})\\ & \ & \ \ \  -  (e^{3r} e^{2 c'_o} - 2 e^r + e^{-r} e^{-2 c'_o} 
+ e^{r}
e^{2 c'_o}  - 2 e^{-r} + e^{-3r} e^{-2 c'_o})\\ & \ & \ \ \ - (e^{3r} - 
e^r -
e^{-r} + e^{-3r})]\\ & \ge & {1\over {8 f^3 (r)}} [C' (C, \eta, c_o) e^{3 
r}] \
\ \ \ \ \ \ \ \
\ \ \ \ \ \ \ \ \ \ \ \ \ \ \ \ \ \ \ \  \ \ \ \ \ \ \ \ \ \ \ \ \ \ \ \ \
\ \ \ \ \ \ 
 \end{eqnarray*} where $C' (C, \eta, c_o)$ is a positive constant. As $f 
(r) =
(e^r + e^{-r})/2$, therefore we can find a point $r_8 > r_7$ such that  
$${{2 g (\alpha (r)) g' (\alpha (r)) \alpha' (r)}\over {f^2 (r)}}  - {{f'
(r)}\over f (r)} [ {{g^2 (\alpha (r))}\over {f^2 (r)}} + (\alpha' (r))^2 
] \ge 
C'' (\eta, c_o) \leqno (4.17)$$
for $r > r_7$ and for some constant $C'' (\eta, c_o) > 0$. That is, 
$$
\Theta' (r) [ ( p - 2 ) (\alpha' (r))^2 + (n - 1) {{g^2 (\alpha 
(r))}\over {f^2
(r)}} ] 
\ge  (p -
2) (n - 1) C'' (\eta, c_o)  \Theta (r)
$$ 
for all $r > r_8$. Since 
$1 - \eta < (\alpha' (r))^2 < 1$ and $C_o \le  g^2 (\alpha (r))/ f^2 (r)) \le
C'_o$, as $r <
\alpha (r) \le r + C$ for all $r > r_8$, we have $C_1 < \theta (r) < 
C'_1$ for
all $r > r_8$ and for some positive constant
$C_o,  C'_o, C_1$ and
$C'_1$. Furthermore $({{g^2 (\alpha)}\over {f^2}})' < 0$ as
in the proof of lemma 4.1. Therefore $\alpha'' (r) \ge c_1$ if $r > r_7$ and
$\alpha' (r) \ge 1 - \eta$, where $c_1$ is a positive constant. For the
sequence $\{x_i\}$, we choose a point $x_i > r_7$ and at $x_i$ we have
$\alpha (x_i) > x_i + c_o$ and $1 > \alpha (x_i) \ge 1 - \eta$. Therefore
$\alpha'' (x_i) \ge c_1$. And for all $r > x_i$ we have $\alpha (x_i) > 
x_i +
c_o$ and $\alpha (x_i) \ge 1 - \eta$ as $\alpha'' (r) \ge c_1 > 0$ for 
all $r
\ge r_i$. Thus $\alpha'' (r) \ge c_1$ for all $r \ge r_i$. Then we can find
$r_8 > r_i$ such that
$\alpha (r_8) = 1$, which is impossible. Hence we must have $c_o = 0$. 
Similar
argument works for $\alpha (r) < r$, using theorem 4.11. \ \ \ \ \ {\bf
Q.E.D.}\\[0.2in]
\hspace*{0.5in}For $p > 2$, let $F (r, \varphi) = (\alpha (r),
\varphi)$ be a rotationally symmetric $p$-harmonic map from the 
hyperbolic space
to itself, where $\alpha \in C^2 (0, \infty)$ and $\lim_{r \to 0^+}
\alpha (r) = \alpha (0) = 0$. By lemma 4.3, there is a point $r_o > 0$ such
that either $\alpha (r) \ge r$ or $\alpha (r) \le r$ for all $r > r_o$. If
$\alpha (r) \not\equiv r$ for $r > r_o$, then there is a point $r_i > r_o$
such that either $\alpha (r_i) > r_i$ or $\alpha (r_i) < r_i$. If $\alpha 
(r_i)
> r_i$ and $\alpha' (r_i) > 1$, then by theorem 4.4, $\alpha (r) > r + C$ for
all $r$ large. If $\alpha (r_i) < r_i$ and $\alpha' (r_1) < 1$, then $\alpha$
is bounded. The other cases are $\alpha' (r) \le 1$
for all $r \ge r_1$ with $\alpha (r) > r$ and $\alpha' (r) \ge 1$ for all 
$r \ge
r_1$ with $\alpha (r) < r$, that is, $\alpha (r)$ is bounded between $r + C$
and $r - C$ for some $C > 0$ and for all $r$ large. Hence $\alpha$ is 
asymptotic
to the line
$y = x$ by theorem 4.15. Therefore we have the following.\\[0.2in]
{\bf Theorem 4.18.} \ \ {\it For $p > 2$, let $F (r, \varphi) = (\alpha (r),
\varphi)$ be a rotationally symmetric $p$-harmonic map from the 
hyperbolic space
to itself, where $\alpha \in C^2 (0, \infty)$ and $\lim_{r \to 0^+}
\alpha (r) = \alpha (0) = 0$. Then either 1) for any constant $C > 0$, 
$\alpha
(r) > r + C$ for all $r$ large, or 2) $\alpha$ is bounded, or 3) $\alpha$ is
asymptotic to the identity map.}

\vspace{0.5in}

{\bf \Large 5. \ \ Asymptotic properties}

\vspace{0.3in}

Let $f, g$ be functions in $C^2 ([0, \infty))$ which satisfy the 
conditions in
(1.1).  We study asymptotic properties of solutions to (1.4). In this 
section we
assume that 
$f$ is similar to the hyperbolic space and $g$ is similar to the Euclidean
space. The result in this section can be founded in [9] We assume that 
there are
constants
$a, C > 0$ such that 
$${1\over C} e^{ar} \le f (r) \le C e^{ar}\,, \ \ \ \ {1\over C} e^{ar} 
\le f' 
(r) \le C e^{ar} \leqno {(\mbox A)}$$ for all $r > 1$. And for $y > 1$, we
assume that  
$$g (y) \le C' y^m\,, \ \ \ \ 0 \le g' (y) \le C' y^{m - 1} \leqno ({\mbox
{B}})$$  for some constants $m \ge 1$ and $C' > 0$.\\[0.2in]  {\bf Lemma 
5.1.} \
\ {\it For $p > 2$, let $\alpha (r) \in C^2 (0, \infty)$ be a positive 
solution
to (1.4) with $\lim_{r \to 0^+} \alpha (r) = 0$. Suppose that
$f$ and
$g$ satisfy the conditions (A) and (B), respectively. For $m > 1$. let $c 
< a/(m
- 1)$ be a positive constant.  Then there exists a positive constant
$r_o$ such that either $\alpha (r) \ge e^{cr}$ or $\alpha (r) \le e^{cr}$ for
all $r > r_o$.}\\[0.1in]  {\bf Proof.} \ \  If $\alpha$ is bounded, then 
we have
$\alpha (r) \le e^{cr}$ for $r$ large enough. Therefore we may assume that
$\lim_{r \to \infty}
\alpha (r) =
\infty$. Hence there exists a constant $r' > 1$ such that $\alpha (r) > 
1$ for
all $r > r'$. For
$m > 1$, suppose that we cannot find such a
$r_o$ as in the statement of the theorem, then for any
$r_o > r'$, there exist two points $r_2 > r_1 > r_o$ such that $\alpha 
(r_1) =
e^{c r_1}$, 
$\alpha (r_2) = e^{c r_2}$ and $\alpha (r') \le e^{c r'}$ for $r' \in [r_1,
r_2]$. There is a point $r \in [r_1, r_2]$ such that $\alpha (r) \le e^{cr}$,
$\alpha' (r) = c e^{cr}$ and $\alpha'' (r) \ge c^2 e^{cr}$. If $c < a/(m 
- 1)$,
we have
\begin{eqnarray*} (5.2) \ \ \ \ \ \ \ \ \ \ \ \ \ \ & \ & {{2 g (\alpha 
(r) ) g'
(\alpha (r))
\alpha' (r)}\over {f^2 (r)}} - {{f' (r)}\over {f (r)}} [ {{ g^2 (\alpha
(r))}\over {f^2 (r)}} + (\alpha' (r))^2 ]\\ & \le & {{2 C'^2 C^2 c e^{2m c
r}}\over { e^{2ar} }}   - {{c^2 e^{2cr} }\over {C^2}} < - c' \ \ \ \ \ \ 
\ \ \ \
\ \ \ \ \ \ \ \ \ \ \ \ \ \ \ \ \ \ \ \ \ \ \
\ \ \ \ \ \ \ 
\end{eqnarray*} if $r_o$ is large enough. Here $c'$ is a positive 
constant. By
(2.3) and (5.2) we have
$$\Theta' (r) [ ( p - 1 ) (\alpha' (r))^2 + (n - 1) {{g^2 (\alpha (r) )}\over
{f^2 (r)}} ] \le - c' (p - 2) (n - 1)
\Theta (r)\,. \leqno (5.3)$$ And 
\begin{eqnarray*} (5.4) \ \ \ \ \ \ \ \ \ \ \ \ \ \ \ \ & \ & \Theta' (r) 
[ ( p
- 1 ) (\alpha' (r))^2 + (n - 1) {{g^2 (\alpha (r) )}\over {f^2 (r)}} ]\\ 
& = &
({p\over 2} - 1) \theta^{{p\over 2} - 2} (r) 
 [ ( p - 1 ) (\alpha' (r))^2 + (n - 1) {{g^2 (\alpha (r) )}\over {f^2 
(r)}} ]\\
& \ & \ \ \ \ \ \ \times [(n - 1) ({{g^2 (\alpha (r))}\over {f^2 (r)}})' 
+ 2
\alpha' (r) \alpha'' (r)]\\ & \ge & ({p\over 2} - 1) \Theta (r) [(n - 1) 
({{g^2
(\alpha (r))}\over {f^2 (r)}})' + 2 \alpha' (r) \alpha'' (r)]\,, \ \ \ \ 
\ \ \ \
\ \ \ \ \ \ \ \ \ \ \ \
\
\end{eqnarray*} as $\theta (r) = (n - 1) {{g^2 (\alpha (r))}\over {f^2 
(r)}} +
(\alpha' (r))^2$ and $\Theta (r) = \theta^{ {p\over 2} - 1} (r)$. By (5.3),
there is a constant $c'' > 0$ such that at $r$ we have
$$(n - 1) ({{g^2 (\alpha (r))}\over {f^2 (r)}})' + 2 \alpha' (r) \alpha'' (r)
\le - c''\,.\leqno (5.5)$$ On the other hand
\begin{eqnarray*} (5.6) \ \ \ \ \ \ \ \ \  \ \ \ \ \ \ \ \ \ \ \ \ \ \ ( 
{{g^2
(\alpha)} \over {f^2 (r)}} )' & = & {{2 f (r) g (\alpha) [ f(r) g' (\alpha)
\alpha' (r) - f' (r) g (\alpha) ]}\over {f^4 (r)}}\\  & \ge & - {{2 f' 
(r) g^2
(\alpha)}\over {f^3 (r)}}\\ & \ge & - 2C'^2  C^2 {{ e^{2mcr} }\over 
{e^{2ar}}}\\
& \ge & - \epsilon  e^{2cr}\,,\ \ \ \ \ \ \ \ \ \ \ \ \ \ \ \ \ \ \ \ \ \ 
\ \ \
\ \ \ \ \ \ \ \ \ \ \ \ \ \ \ \ \ \ \ \ \ \ \ \
\end{eqnarray*}   as $c < a/(m - 1)$. Here $\epsilon$ is a positive constant.
And by choosing
$r_o$ large, we may assume that $\epsilon < 2c^3/(n - 1)$. There if
$r_o$ is large and
$c < a/(m - 1)$, then (5.5) and (5.6) imply that 
$$2 \alpha' (r) \alpha'' (r) \le  (n - 1) \epsilon c^{2cr}\,.$$ That is 
$$\alpha'' (r) \le {{(n - 1)\epsilon} \over {2c}} c^{cr} < c^2 
e^{cr}\,,$$ which
is a contradiction. Hence there exists a positive constant
$r_o > r'$ such that either $\alpha (r) \ge e^{cr}$ or $\alpha (r) \le 
e^{cr}$
for all $r > r_o$.\ \ \ \ \ {\bf Q.E.D.}\\[0.2in]  {\bf Lemma 5.7.} \ \ 
{\it For
$p > 2$, let $\alpha (r) \in C^2 (0, \infty)$ be a positive solution to (1.4)
with $\lim_{r \to 0^+} \alpha (r) = 0$. Suppose that
$f$ and
$g$ satisfy the conditions (A) and (B) with $m > 1$, respectively. Let $c 
< a/(m
- 1)$ be a positive constant. Suppose that there exists a $r_o > 0$ such that
$\alpha (r)
\le e^{cr}$ for all $r > r_o$, then $\alpha (r) \le C$ for some positive 
constant
$C$.}\\[0.1in] {\bf Proof.} \ \  Assume that $\alpha (r)
\le e^{cr}$ for all
$r > r_o$, where
$c < a/(m - 1)$. We may assume that $\alpha \not\equiv 0$. We can find 
positive 
constants $r'$ and $\delta$ such that $\alpha (r) \ge \delta $ for all $r
\ge r'$. There is a positive constant $C''$ such that 
$$g (y) \le C'' y^m\,, \ \ \ \ 0 \le g' (y) \le C'' y^{m - 1}$$ for all 
$y \ge
\delta$. We may take $r_o > r' > 1$. Then for
$r > r_o$ and
$1 >
\tau > 0$, we have
\begin{eqnarray*} (5.8) \ \ & \ & {{2 g (\alpha (r)) g' (\alpha (r))
\alpha' (r)}\over {f^2 (r)}} - {{f' (r)}\over {f (r)}} [ {{ g^2 (\alpha
(r))}\over {f^2 (r)}} + (\alpha' (r))^2 ]\\ & = & {{2 g (\alpha (r)) g' 
(\alpha
(r))}\over {f^{3\over 2} (r)f'^{1\over 2} (r)}} [ {{f'(r) (\alpha' 
(r))^2}\over
{f (r)}} ]^{1\over 2} - {{ f' (r) g^2 (\alpha (r))}\over {f^3 (r)}} - {{f'
(r)}\over {f (r)}} (\alpha' (r))^2\\  & \le & {1\over {1 - \tau}} {{g^2 
(\alpha
(r)) g'^2 (\alpha (r))}\over {f^3 (r) f' (r)}} + (1 - \tau) {{f' 
(r)}\over {f
(r)}} (\alpha' (r))^2 - {{ f' (r) g^2 (\alpha (r))}\over {f^3 (r)}} - {{f'
(r)}\over {f (r)}} (\alpha' (r))^2\\  & = & {{g^2 (\alpha (r)) }\over {(1 -
\tau) f^3 (r) f' (r)}} [g'^2 (\alpha (r)) - ( 1-
\tau)^2 f' (r)^2] - \tau [{{ f' (r) g^2 (\alpha (r))}\over {f^3 (r)}} + {{f'
(r)}\over {f (r)}} (\alpha' (r))^2]\\  & \le & {{g^2 (\alpha (r)) }\over 
{(1 -
\tau) f^3 (r) f' (r)}}  [C''^2 e^{2 c (m - 1) r} - {{( 1 - \tau )^2}\over 
{C^2}}
e^{2ar}] - \tau [{{ f' (r) g^2 (\alpha (r))}\over {f^3 (r)}} + {{f' 
(r)}\over {f
(r)}} (\alpha' (r))^2]\,, 
\end{eqnarray*} where we have used the inequality $2AB \le (1 - \tau) A^2 
+ 1/(1
- \tau) B^2$ for $0 < \tau < 1$. Hence there exists a constant $r'' > 
r_o$ such
that for $r > r''$ we have 
$$C''^2 e^{2 c (m - 1) r} - {{( 1 - \tau )^2}\over {C^2}} e^{2ar} \le 
0\,.$$ For
$r > r''$ we have
\begin{eqnarray*} (5.9) \ \ \ \ \ \ \ \ \ \ \ \ \ \ \  & \ & {{2 g 
(\alpha (r))
g' (\alpha (r))
\alpha' (r)}\over {f^2 (r)}} - {{f' (r)}\over {f (r)}} [ {{ g^2 (\alpha
(r))}\over {f^2 (r)}} + (\alpha' (r))^2 ]\\ & = & - \tau {{f' (r)}\over 
{f (r)}}
[ {{g^2 (\alpha (r))}\over {f^2 (r)}} + (\alpha' (r))^2]\\ & \le & 
{{-\tau}\over
{C^2}}[ {{g^2 (\alpha (r))}\over {f^2 (r)}} + (\alpha' (r))^2]\,. \ \ \ \ 
\ \ \
\ \ \ \ \ \ \ \ \ \ \ \ \ \ \ \ \ \ \ \ \ \ \ \ \ \ \ \
\ \ \ \ \ \ \ \ \ \  
\end{eqnarray*} Together with (2.3) we have
\begin{eqnarray*} (5.10) \ \ \ \ \ \ \ \ \ \ \ \ \ \ \ & \ & \Theta' (r) 
[ (p -
1) (\alpha' (r))^2 + (n - 1) {{g^2 (\alpha)}\over {f^2 (r)}}]\\ & \le & -
{{\tau}\over {C^2}} (p - 2) (n - 1) \, [ {{g^2 (\alpha (r))}\over {f^2 
(r)}} +
(\alpha' (r))^2] \Theta (r)\,.  \ \ \ \ \ \ \
 \ \ \ \ \ \ \  \ \ \ \ \ \ \  \ \ \ \ \ \ \
\end{eqnarray*} Thus we have 
$$\Theta' (r) \le - c' \Theta (r) \ \ \ \ {\mbox {for \ all}} \ \ \ r >
 r''\,, \leqno (5.11)$$  where $c' = (p - 2) (n - 1)\tau/(bC^2)$ and $b = 
\max \
\{ p - 2, n - 1 \}$ are positive constants. Integrating both sides of the 
above
inequality we have
$$\Theta (r) \le e^{-c' (t - \bar r)} \Theta (r'')\,, \ \ \ r > r''\,.$$  
So we
have
$$\alpha' (r) \le C_1 e^{- {{c' (t - r'')}\over {p - 2}}} \ \ \ \ {\mbox 
{for \
all}} \ \ \ \ r > r''\,,$$  where $C_1 = (\Theta ( r''))^{1\over {p - 
2}}$ is a
positive constant. As $\alpha' (r) > 0$, upon integration we conclude that
$\alpha$ is bounded. \ \ \ \ \ {\bf Q.E.D.}\\[0.2in]
\hspace*{0.5in}Combining the above two lemmas we have the following. 
\\[0.2in]
{\bf Theorem 5.12.} \ \ {\it For $p > 2$, let $\alpha (r) \in C^2 (0, 
\infty)$ be
a positive solution to (1.4) with $\lim_{r \to 0^+} \alpha (r) = 0$. 
Suppose that
$f$ and
$g$ satisfy the conditions (A) and (B) with $m > 1$, respectively. Let $c 
< a/(m
- 1)$ be a positive constant. Then either there exists a positive 
constant $r_o$
such that
$\alpha (r) \ge e^{cr}$ for all $r > r_o$ or there is a finite positive 
number
$C$ such that $\alpha (r) \le C$ for all $r > 0$.}\\[0.2in] 
\hspace{0.5in}In the following we assume that  
$$ g' (y) \le C_2 g (y) \ \ \ \ \ \ {\mbox { for \ all }} \ \ y \ge 1\,, 
\leqno
({\mbox {C}})
$$ where $C_2$ is a positive constant. Condition (C) implies that $g$ 
grows at
most exponentially. Condition (B) is stronger than condition (C). The 
following
could be considered as a generalization of lemma 5.7.\\[0.2in] 
{\bf Lemma 5.13.}
\ \ {\it For
$p > 2$, let
$\alpha (r)
\in C^2 (0,
\infty)$ be a positive solution to (1.4) with $\lim_{r \to 0^+} \alpha 
(r) = 0$.
Suppose that $f$ and $g$ satisfy the conditions (A) and (C), 
respectively, and
$g' (y) > 0$ for all $y > 0$. There is a positive constant $c > 0$ such 
that if
there exists a positive constant
$r_o$ such that $\alpha (r)
\le g^{-1} ( c f (r) )$ for all
$r > r_o$, then $\alpha$ is a bounded function on ${\real}^+$.}\\[0.2in] {\bf
Proof.} \ \ We can find positive  constants $r'$ and $\delta$ such that
$\alpha (r) \ge \delta $ for all $r \ge r'$. There is a positive constant 
$C_3$
such that
$$  g' (y) \le C_3 g (y) \ \ \ \ \ \ {\mbox { for \ all }} \ \ y \ge 
\delta\,.
$$ For $0 < \tau < 1$ and $r > r'$, from (5.8) we have 
\begin{eqnarray*} & \ & \ {{2 g (\alpha (r)) g' (\alpha (r))
\alpha' (r)}\over {f^2 (r)}} - {{f' (r)}\over {f (r)}} [ {{ g^2 (\alpha
(r))}\over {f^2 (r)}} + (\alpha' (r))^2 ]\\  & \le & {{g^2 (\alpha (r)) 
}\over
{(1 -
\tau) f^3 (r) f' (r)}} [g'^2 (\alpha (r)) - ( 1-
\tau)^2 f' (r)^2] - \tau [{{ f' (r) g^2 (\alpha (r))}\over {f^3 (r)}} + {{f'
(r)}\over {f (r)}} (\alpha' (r))^2]\\  & \le & {{g^2 (\alpha (r)) }\over 
{(1 -
\tau) f^3 (r) f' (r)}} [C_3^2 g^2 (\alpha (r)) - {{(1 - \tau)^2}\over 
{C^2}} f^2
(r)] - \tau [{{ f' (r) g^2 (\alpha (r))}\over {f^3 (r)}} + {{f' (r)}\over {f
(r)}} (\alpha' (r))^2]\,.
\end{eqnarray*} If we choose a positive constant $c$ such that $\,c \le 
(1 -
\tau) /(CC_3)\,,$ and if $\alpha (r) \le g^{-1} ( c f (r) ),$ then we 
have 
$${{2 g (\alpha (r)) g' (\alpha (r))
\alpha' (r)}\over {f^2 (r)}} - {{f' (r)}\over {f (r)}} [ {{ g^2 (\alpha
(r))}\over {f^2 (r)}} + (\alpha' (r))^2 ] \le  - \tau [{{ f' (r) g^2 (\alpha
(r))}\over {f^3 (r)}} + {{f' (r)}\over {f (r)}} (\alpha' (r))^2]\,.$$ We can
process as in the proof of lemma 5.7.\qed
\hspace*{0.5in}The above lemma provides asymptotic information about
$\alpha (r)$. We can also obtain medium range information about
$\alpha (r)$.\\[0.2in]  {\bf Lemma 5.14} \ \ {\it For
$p > 2$, let
$\alpha (r)
\in C^2 (0,
\infty)$ be a positive solution to (1.4) with $\lim_{r \to 0^+} \alpha 
(r) = 0$.
Suppose that
$f$ and
$g$ satisfy the conditions (A) and (C), respectively, and $g' (y) \ge 0$ 
for all
$y > 0$. Assume that
$\alpha$ is not a bounded function on ${\bf R}^+$ and let $\alpha (\bar 
r) = 1$
for some
$\bar r > 0$. Then there exists a positive constant
$\delta$ such that the energy density $\theta (r) \ge \delta$ for all $r 
\ge \bar
r$.}
\\[0.2in] {\bf Proof.} \ \ From (2.3) we have
\begin{eqnarray*} (5.15)& \ & \Theta' (r) [ ( p - 1 ) (\alpha' (r))^2 + 
(n - 1)
{{g^2 (\alpha (r) )}\over {f^2 (r)}} ]\\ & = & (p - 2) 
\Theta (r) \{ {{2 (n - 1)g (\alpha (r) ) g' (\alpha (r))
\alpha' (r)}\over {f^2 (r)}} - {1\over 2} {{f' (r)}\over {f (r)}} [ {{ g^2
(\alpha (r))}\over {f^2 (r)}} + (\alpha' (r))^2]\\ & \ & \ \ \ \  - (n - 
{3\over
2}){{f' (r)}\over {f (r)}} [ {{ g^2 (\alpha (r))}\over {f^2 (r)}} + (\alpha'
(r))^2]
\}\,.
\end{eqnarray*} At a point $r \ge \bar r$, we have $\alpha (r) \ge 1$. So 
at the
point $r > \bar r$, we have either
$${{2 (n - 1)g (\alpha (r) ) g' (\alpha (r))
\alpha' (r)}\over {f^2 (r)}} \ge {1\over 2} {{f' (r)}\over {f (r)}} [ {{ g^2
(\alpha (r))}\over {f^2 (r)}} + (\alpha' (r))^2]\,, \leqno ({\mbox i})$$ or
\begin{eqnarray*} ({\mbox {ii}}) \ \ \ \ \ \ \ \ \ & \ & \Theta' (r) [ ( 
p - 1 )
(\alpha' (r))^2 + (n - 1) {{g^2 (\alpha (r) )}\over {f^2 (r)}} ]\\ & \le 
& - (p
- 2) (n - {3\over 2}) \Theta (r) {{f' (r)}\over {f (r)}} [ {{ g^2 (\alpha
(r))}\over {f^2 (r)}} + (\alpha' (r))^2]\,.  \ \ \ \ \ \ \ \ \ 
\ \ \ \ \ \ \ \ \  \ \ \ \ \ \ \ \ \  \ \ \ \ \ \ \ \ \  \ \ \ \ \ \ \ \
\  \ \ \ \ \ \ \ \ \  \ \ \ \ \ \ \ \ \ 
\end{eqnarray*} In case of (i), from condition (A) and (C) we have
$$2C_2 (n - 1) \alpha' (r) \ge {{f' (r) f (r) (\alpha' (r))^2}\over {g^2 
(\alpha
(r))}} + {{f' (r)}\over {2 f (r)}} \ge {1\over {2 C^2}}\,.$$ By choosing a
bigger constant if necessary, we may assume that (A) holds for
$r > \bar r$. Thus if (i) holds at
$r
\ge
\bar r$, we have
$$\alpha' (r) \ge {1\over {4 (n - 1) C_2C^2}}\,. \leqno (5.16)$$ Take 
$$\delta = ({1\over {8(n - 1) C_2C^2 }})^2\,. \leqno (5.17)$$ Suppose 
that there
is a point $r'' \ge r_o$ such that $\theta (r'') \le \delta$.  Let 
$$O = \{ t \in [r'', \infty) \  | \ \theta (r) \le \delta \ \ {\mbox {for 
\ all}}
\ \ r \in [r'', t] \}\,.$$ Therefore $r' \in 0$ and
$$
\alpha' (r') \le {1\over {8 (n - 1) C_2C^2}}\,.
$$ We can find $t$ such that  $t > r'$ and  
$$\alpha' (r) \le {1\over {4 (n - 1) C''C^2}}$$ for all $r \in [r'', t]\,.$
Hence we have (ii)  for all $r
\in [r'', t]\,.$ Thus $\Theta' (r) \le 0$ and hence $\theta' (r) \le 0$ 
for all
$r \in [r'', t]\,.$ We have $\theta (r) \le \delta$ for all $r \in [r'', 
t]\,.$
Suppose that sup $O = a <
\infty$. Then at the point $a$ we have
$\theta (a) \le \delta$. The same argument shows that there exists a positive
number $\epsilon$ such that $\Theta' (r) \le 0$ for all $r \in (a, a +
\epsilon)$. Then $a + \epsilon \in O$, which is a contradiction. 
Therefore we
have sup $O =  \infty$, or    
$$\theta (r) \le \delta \ \ {\mbox {for \ all}} \ \ r \in [r', \infty)\,.
\leqno (5.18)$$ In particular, 
$$\alpha' (r) < {1\over {4 (n - 1) C''C^2}}$$ for all $r \in [r', 
\infty)$. By
(ii) we have
$$\Theta' (r) [ ( p - 1 ) (\alpha' (r))^2 + (n - 1) {{g^2 (\alpha (r) )}\over
{f^2 (r)}} ] \le - {{(p - 2) (n - {3\over 2})}\over {C^2}} \Theta [ {{ g^2
(\alpha (r))}\over {f^2 (r)}} + (\alpha' (r))^2]\,.$$ A similar argument 
as in
the proof of lemma 5.7 shows that there is a positive constant $c > 0$ 
such that
$\alpha (r) \le c$ for all $r \in (0,
\infty)$. This contradicts the assumption that $\alpha$ is not  bounded.
Therefore $\theta (r) \ge \delta$ for all $r \ge \bar r$. \ \ \ \ \ {\bf
Q.E.D.}\\[0.2in] 
{\bf Theorem 5.19} \ \ {\it For $p > 2$, let $\alpha (r)
\in C^2 (0,
\infty)$ be a positive solution to (1.4) with $\lim_{r \to 0^+} \alpha 
(r) = 0$.
Suppose that
$f$ and
$g$ satisfy the conditions (A) and (C), respectively, and $g' (y) > 0$ 
for all
$y > 0$. Then either there exist positive constants $r_o$ and $c_o$ such that
either
$\alpha (r) \ge c_o (r - r_o)$ for all
$r > r_o\,,$ or $\alpha$ is a bounded function on ${\bf R}^+$.}\\[0.2in] {\bf
Proof.} \ \ If $\alpha$ is not bounded, then there is a positive constant 
$\bar
r$ such that $\alpha (\bar r) = 1$. By lemma 5.14, we have $\theta (r) \ge
\delta$ for all
$r \ge \bar r$. With the conditions (A) and (C), we can find positive 
constants
$c'$ and $r'_o$ such that 
$$g^{-1} ({\delta \over {\sqrt {2 (n - 1)} }}  f (r)) \ge c' r$$ for all 
$r >
r'_o$. We can choose $c_o < \min \{ c', \sqrt {\delta / 2} \}$ and
$r_o > \max \{ \bar r\,, r'_o \}$. If there is a point $r' > r_o$ such that
$\alpha (r') < c_o(r' - r_o)$. Then we can find a point $r > r_o$ such that
$\alpha (r) < c_o r$ and
$\alpha' (r) < c_o$. Thus
$$0 < \alpha'(r') < \sqrt {\delta \over 2}\,.$$ And 
$$ (n - 1) {{g^2 (\alpha (r'))}\over {f^2 (r)}} \le  (n - 1) {{g^2 (cr)}\over
{f^2 (r)}} \le (n - 1) {{g^2 (g^{-1} ({\delta \over {\sqrt {2 (n - 1)} }} 
f (r))
)}\over {f^2 (r)}} \le {\delta \over 2}\,,
$$ as $g' \ge 0$ implies that $g$ is non-decreasing. Therefore $\theta < 
\delta$,
contradiction.\ \ \ \ \ {\bf Q.E.D.}\\[0.2in]  

\vspace{0.5in}

{\bf \Large 6. \ \ Metrics with polynomial growth}

\vspace{0.3in}

For $p > 2$, let $F (r, \varphi) = (\alpha (r),
\varphi)$ be a rotationally symmetric $p$-harmonic map from $M^n (f)$ to
itself. The results in this section can be founded in [8].\\[0.2in] 
{\bf Lemma
6.1.}
\
\ {\it For
$f (r) = r^m$ with
$m
\ge 1$, if there exists a point $r_o > 0$ such that $\alpha (r_o) < cr$ and
$\alpha' (r_o) < c$ for some $c \in (0, 1]$, then $\alpha (r) < cr$ for 
all $r
\ge r_o$.}\\[0.2in] {\bf Proof.} \ \ Assume that there is a point $r' > 
r_o$ such
that $\alpha (r') = cr'$. As $\alpha' (r_o) < c$, there exists a point
$r \in (r_o, r_1]$ such that
$\alpha (r) = c'r$ for some constant $c' \in (0, c)$ and $\alpha' (r) = 
c$ and
$\alpha'' (r)
\ge 0$. At
$r$, we have    
\begin{eqnarray*} 
(6.2) \ \ \ \ \ \ \ \ \ \ \ \ \ \ \ \ \ & \ & \ {{2 g (\alpha (r)) g' (\alpha
(r))
\alpha' (r)}\over {f^2 (r)}} - {{f' (r)}\over {f (r)}} [ {{ g^2 (\alpha
(r))}\over {f^2 (r)}} + (\alpha' (r))^2 ]\\ & = & {{2m \alpha^{2m - 1}c}\over
{r^{2m}}} - {m\over r} ({{\alpha^{2m}}\over {r^{2m}}} + c^2)\\ & = & 
{m\over r}
(2 c c'^{2m - 1} - c'^{2m} - c^2) \ \ \ \ \ \ \ \ \ \ \ \
\ \ \ \ \ \ \ \ \ \ \ \ \ \ \ \ \ \ \ \ \ \ \ \ \ \ \ \ \ \ 
\end{eqnarray*} It can be shown that if $m \ge 1$, then $2 c c'^{2m - 1} -
c'^{2m} - c^2 < 0$ for
$c'
\in (0, c)$. In fact, if we let 
$$\phi (x) = 2 c x^{2m - 1} - x^{2m} - c^2\,,$$ then $\phi (c) = c^{2m} - c^2
\le 0$ as $c \le 1$, $m \ge 1$, and $\phi' (x) > 0$ for $x \in (0, c)$. Hence
$\Theta' (r) < 0$. Therefore 
$$0 > \theta' (r) = ( {{g^2 (\alpha)}\over {f^2 (r) }})' + 2 \alpha' (r) 
\alpha''
(r)\,.$$  And
$$({{g^2 (\alpha)}\over {f^2 (r)}})' = {{2m}\over {r^{4m}}} (r^{2m} 
\alpha^{2m -
1} c - \alpha^{2m} r^{2m - 1}) = {{2m}\over r} (c'^{2m - 1} c - c'^{2m}) 
> 0
\leqno (6.3)$$  as $0 < c' < c$. Hence we have $\alpha'' (r) > 0$, 
contradiction.
\ \ \ \ \ {\bf Q.E.D.} \\[0.2in]  {\bf Theorem 6.4.} \ \ {\it  Let $f (r) 
= g
(r) = r^m$ for some
$m
\ge 1$. Let
$\alpha \in C^1 [0, \infty) \cap C^2 (0, \infty)$ be a solution to the
rotationally symmetric $p$-harmonic equation with $\alpha (0) = 0$. If 
$\alpha'
(0) = c$ for some $c \in (0, 1]$, then
$\alpha (r) \le cr$ for all $r > 0$.}\\[0.2in] {\bf Proof.} \ \ For 
$\alpha' (0)
= c < 1$, assume that there is a point $r' > 0$ such that $\alpha (r') > 
c r'$.
Then we can find $c' > c$ such that $\alpha (r') > c' r'$. Since $\alpha' 
(0) =
c < c'$. we can find a point $r_o < r'$ close to zero such that $\alpha 
(r_o) <
c' r_o$ and $\alpha' (r_o) < c'$. Lemma 6.1 implies that $\alpha (r) <  
c'r$ for
all $r > r_o$, contradiction. So we have $\alpha (r) \le cr$ for all $r > 
0$. In
case $c = 1$, $\alpha (r) = r$ is the unique solution to the equation with
$\alpha' (0) = 1$ [2]. \ \ \ \ \ {\bf Q.E.D.}\\[0.2in]
\hspace*{0.5in}For $m = 1$, that is, the Euclidean space, and for any 
positive
number $c$, the function
$\alpha_c = cr$ is a rotationally symmetric $p$-harmonic map into the 
Euclidean
space for all $p > 0$.\\[0.2in] {\bf Theorem 6.5.} \ \ {\it Let $f (r) = r^m$
with $m > 1$. For $p > 2$, let $F (r, \varphi) = (\alpha (r),
\varphi)$ be a rotationally symmetric $p$-harmonic map from $M^n (f)$ to 
itself,
where $\alpha
\in C^2 (0, \infty)$ and $\lim_{r \to 0^+}
\alpha (r) = \alpha (0) = 0$. Then for any $c > 0$, we can find a positive
number $r_c$ such that either $\alpha (r) \ge c r$ or $\alpha (r) \le cr$ for
all $r > r_c$.}\\[0.2in] {\bf Proof.} \ \ Let  $c \in (0, 1]$ and suppose 
that
there exist two points $r_2 > r_1 > 0$ such that $\alpha (r_1) = cr_1$, 
$\alpha
(r_2) = cr_2$ and $\alpha (r) \le cr$ for all $r \in [r_1, r_2]$ but 
$\alpha (r)
\not\equiv cr$ on $[r_1, r_2]$. Then there exist a point $r_o \in (r_1, r_2)$
such that $\alpha (r_o) < cr_o,
\alpha' (r_o) < c$. Lemma 6.1 implies $\alpha (r) < cr$ for all $r \ge r_o$,
which is a contradiction. Hence we can find a positive number $r_c$ such that
either $\alpha (r) \ge c r$ or $\alpha (r) \le cr$ for all $r > r_c$. Let 
$c > 1$ and suppose that there exist two points $r_2 > r_1 > 0$ such that
$\alpha (r_1) = cr_1$, $\alpha (r_2) = cr_2$ and $\alpha (r) \ge cr$ for 
all $r
\in [r_1, r_2]$ but $\alpha (r) \not\equiv cr$ on $[r_1, r_2]$. Let 
$$c_o = \inf \ \{ \gamma \ge c \ | \ \gamma r \ge \alpha (r) \ \ {\mbox {for
\ all}} \ \ r \in [r_1, r_2] \}\,.$$ Thus $c_o r \ge \alpha (r)$ for all  
$r \in
[r_1, r_2]$. Let 
$$c' = {\mbox {max}} \ \{c\,, {m \over {2m - 1}} c_o \}\,.$$  As $m > 1$ and
$c_o > c$, we have
$c' < c_o$. Hence there exist $x'_1, x'_2$ such that
$x_1 \le x'_1 < x'_2 \le x_2$ and $\alpha (r_1) = c'r_1$, $\alpha (r_2) = 
c'r_2$
and $\alpha (r) \ge c'r$ for all $r \in [r'_1, r'_2]$. Then there exists 
a point
$r_o \in [r'_1, r'_2]$ such that $\alpha (r_o) \ge c'r_o, \alpha' (r_o) = c'$
and $\alpha'' (r_o) \le 0$. Let $\alpha (r_o) = c'' r_o$ for some
$c'' \ge c'$. As in (6.2) we have
\begin{eqnarray*} & \ & \ {{2 g (\alpha (r_o)) g' (\alpha (r_o)) \alpha'
(r_o)}\over {f^2 (r_o)}} - {{f'(r_o)}\over {f(r_o)}} [ {{ g^2 (\alpha
(r_o))}\over {f^2 ((r_o))}} + (\alpha' (r_o))^2 ]\\ & = & {m\over t} (2 c'
c''^{2m - 1} - c''^{2m} - c'^2)\,.
\end{eqnarray*}  The function $\phi (x) = 2 c' x^{2m - 1} - x^{2m} - 
c'^2$ has
$\phi (c') = c'^{2m} - c'^2 > 0$, since $c' \ge c > 1$ and $m > 1$. While 
$$\phi' (x) = 2m x^{2m - 2} ({{2m - 1}\over m} c' - x) = {{2m x^{2m - 
2}}\over
r_o} ({{2m - 1}\over m} c'r_o  - xr_o)\,.$$   As $c' \ge m/(2m - 1) c_o$,
therefore 
$${{2m - 1}\over m} c'r \ge c_o r \ge \alpha (r) \ \ \ \ {\mbox {for \ 
all}} \
\ r \in [r'_1, r'_2]\,.$$ In particular 
$${{2m - 1}\over m} c'r_o \ge \alpha (r_o) = c'' r_o \ge x r_o \ \ \ 
{\mbox {for
 \ all}} \ \ x \in [c', c'']\,.$$ Thus $\phi' (x) \ge 0$ for all $x \in [c',
c'']$. Hence
$\phi (c'') > 0$. Therefore 
$${{2 g (\alpha (r_o)) g' (\alpha (r_o)) \alpha' (r_o)}\over {f^2 (r_o)}} -
{{f'(r_o)}\over {f(r_o)}} [ {{ g^2 (\alpha (r_o))}\over {f^2 ((r_o))}} +
(\alpha' (r_o))^2 ] > 0\,. \leqno (6.6)$$  As in (6.3) 
$$({{g^2 (\alpha (r_o))}\over {f^2 (r_o)}})' = {{2m}\over {r_o^{4m}}} 
(r_o^{2m}
\alpha^{2m - 1} (r_o) c' - \alpha^{2m} (r_o) r_o^{2m - 1}) \le 0\,, \leqno
(6.7)$$  as $\alpha (r_o) \ge c' r_o$. Hence $\alpha'' (r_o) > 0$. This 
is a
contradiction. So we can find a positive number $r_c$ such that either 
$\alpha
(r) \ge c r$ or $\alpha (r) \le cr$ for all $r > r_c$. \ \ \ \ \ {\bf
Q.E.D.}\\[0.4in] 
ACKNOWLEDGEMENT. \ \ We would like to thank Professor
J.B. McLeod for valuable discussions.

\pagebreak


\begin{thebibliography}{9}
\bibitem{Ch-Law} Cheung, L.-F., \& Law, C.-K., {\it An intial value 
approach to
rotationally symmetric harmonic maps,} Preprint, 1995.
\bibitem{Ch-Law-Leung}  Cheung, L.-F., Law, C.-K. \& Leung, M.-C., {\it 
Entire
solutions of quasilinear differential equations cooresponding to $p$-harmonic
maps,} Preprint, 1995.
\bibitem{CL} Coddington, E.A. \& Levinson, N., {\it ``Theory of Ordinary
Differential Equations''}, (1955) McGraw-Hill.
\bibitem{Eells-L} Eells, J. \& Lemaire, L., {\it Another report of harmonic
maps,} Bulletin of the London Math. Soc. {\bf 20} (1988), 385-524.
\bibitem{G-W} Greene, R. \& Wu, H., {\it Function Theory on Manifolds Which
Possess a Pole,} Lectures Notes in Math. vol. 699, Springer-Verlag, New York,
1979.
\bibitem{H-K-M} Heinonen, J., Kilpel\"ainen, T. \& Martio, O., {\it Nonlinear
Potential Theory of Degenerate Elliptic Equations,} Oxford Mathematical
Monographs, Oxford University Press, Oxford, 1994. 
\bibitem{L} Leung, M.-C., {\it On the infinitesmal rigidity of harmonic 
maps,}
Preprint, 1995.
\bibitem{L} Leung, M.-C., {\it Asymptotic behavior of rotationally symmetric
$p$-harmonic maps,} to appear.
\bibitem{L} Leung, M.-C., {\it Positive solutions of second order quasilinear
equations corresponding to $p$-harmonic maps,} Preprint, 1996.
\bibitem{N-T}Nakauchi, N. \& Takakuwa, S., {\it   A Remark on p-harmonic 
maps,} 
Nonlinear Analysis, Theory, Method \& Applications {\bf 25} (1995), 169-185.
\bibitem{R-R} Ratto, A. \& Rigoli, M., {\it On the asymptotic behaviour
of rotationally symmetric harmonic maps}, Journal of Diff. Eqn. {\bf 101} 
(1993)
15-27.
\bibitem{Ta} Takakuwa, S., {\it Stability and Liouville theorems of 
p-harmonic
maps}, Japan J. Math., {\bf 17} (1991), 317-332.
\bibitem{U} Uhlenbeck, K., {\it Regularity for a class of nonlinear elliptic
systems}, Acta Math.\ {\bf 138} (1970) 219-240.
\bibitem{Wh} White, B., {\it Homotopy classes in Sobolev spaces and the 
existence
of energy minimizing maps,} Acta Math. {\bf 160} (1988), 1-17. 
\bibitem{X-Y} Xu, X.-W. \& Yang, P., {\it A construction of 
$\,m$-harmonic maps
of spheres,} International J. of Math. {\bf 4} (1993), 521-533.
\bibitem{Y} Yau, S.-T.,  {\it Harmonic functions on complete Riemannian
manifolds,} Comm. Pure Appl. Math. \ {\bf 28} (1975), 201-228.
\end{thebibliography}
\end{document}